


\input amssym  

\newdimen\normalparindent
\newdimen\chapterindent
\newdimen\chapterbaselineskip
\newdimen\sectionindent
\newdimen\sectionbaselineskip
\newdimen\subsectionindent
\newdimen\paperwidth
\newdimen\paperheight


\hsize=15truecm
\hoffset=.46truecm
\vsize=23.7truecm
\voffset=.46truecm

\font\eleveni=cmmi11
\font\eighti=cmmi8
\font\sixi=cmmi6
\font\elevensy=cmsy11
\font\eightsy=cmsy8
\font\sixsy=cmsy6

\font\sevenrm=cmr7
\font\twelvei=cmmi12
\font\ninei=cmmi9
\font\seveni=cmmi7
\font\twelvesy=cmsy12
\font\ninesy=cmsy9
\font\sevensy=cmsy7

\font\sevenbf=cmbx7

\skewchar\eleveni='177 \skewchar\eighti='177 \skewchar\sixi='177
\skewchar\elevensy='60 \skewchar\eightsy='60 \skewchar\sixsy='60
\skewchar\twelvei='177 \skewchar\ninei='177 \skewchar\seveni='177
\skewchar\twelvesy='60 \skewchar\ninesy='60 \skewchar\sevensy='60

\font\fivecal=cmsy5

\font\sevencal=cmsy7

\font\tencal=cmsy10

\newfam\calfam
\newfam\frakfam

\def\tenpoint{%
\textfont0=\tenrm \scriptfont0=\sevenrm \scriptscriptfont0=\fiverm
\def\rm{\fam0\tenrm}%
\textfont1=\teni \scriptfont1=\seveni \scriptscriptfont1=\fivei
\def\mit{\fam1}\def\oldstyle{\fam1\teni}%
\textfont2=\tensy \scriptfont2=\sevensy \scriptscriptfont2=\fivesy
\textfont3=\tenex \scriptfont3=\tenex \scriptscriptfont3=\tenex
\textfont\itfam=\tenit
\def\it{\fam\itfam\tenit}%
\textfont\slfam=\tensl
\def\sl{\fam\slfam\tensl}%
\textfont\bffam=\tenbf \scriptfont\bffam=\sevenbf
\scriptscriptfont\bffam=\fivebf
\def\bf{\fam\bffam\tenbf}%
\textfont\ttfam=\tentt
\def\tt{\fam\ttfam\tentt}%
\font\sc=cmcsc10
\font\chapterfont=cmssbx12 scaled \magstep4
\font\titlefont=cmr12 scaled \magstep3
\font\sectionfont=cmbx10 scaled \magstep2
\font\subsectionfont=cmbx10
\font\cmcyr=cmcyr10
\font\cmcti=cmcti10
\font\cmccsc=cmccsc10
\textfont\calfam=\tencal \scriptfont\calfam=\sevencal
  \scriptscriptfont\calfam=\fivecal
\def\cal{\fam\calfam\tencal}%
\normalparindent=20pt
\chapterindent=30pt
\chapterbaselineskip=30pt
\sectionindent=30pt
\sectionbaselineskip=18pt
\subsectionindent=30pt
\smallskipamount=3pt plus 1pt minus 1pt
\medskipamount=6pt plus 2pt minus 2pt
\bigskipamount=12pt plus 2pt minus 2pt
\normalbaselineskip=12.4pt
\normallineskip=1pt
\normallineskiplimit=0.5pt
\jot=3pt
}


\def\S{\mathhexbox278\thinspace}
\def\SS{\mathhexbox278\mathhexbox278\thinspace}

\def\no{$\rm n^o$\thinspace}
\def\nos{$\rm n^{os}$\ }

\def\square{\hbox to.77778em{%
\hfil\vrule\vbox to.675em{\hrule width.6em\vfil\hrule}\vrule\hfil}}

\def\definition#1\par{\medbreak\noindent{\bf Definition.}\enspace
  #1\par\medbreak}
\def\example#1\par{\medbreak\noindent{\bf Example.}\enspace
  #1\par\medbreak}
\def\acknowledgements#1\par{\medbreak\noindent{\it Acknowledgements\/}.\enspace#1\par\medbreak}
\long\def\remark#1\par{\medbreak\noindent{\it Remark\/}.\enspace#1\par\medbreak}
\long\def\remarks#1\par{\medbreak\noindent{\it Remarks\/}.\enspace#1\par\medbreak}
\def\exercise#1\par{\medbreak\noindent{\bf Exercise.}\enspace
#1\par\medbreak}
\def\notation#1\par{\medbreak\noindent{\bf Notation.}\enspace
#1\par\medbreak}
\def\proof{\noindent{\it Proof\/}.\enspace}
\def\endproof{\nobreak\hfill\quad\square\par\medbreak}

\def\lineover#1{{\offinterlineskip\mathchoice
{\setbox0=\hbox{$\displaystyle#1$}%
\vbox{\kern .33pt\hbox to\wd0{\kern 1pt\leaders\hrule height .33pt%
\hfill\kern 1pt}\kern 1pt\box0}}
{\setbox0=\hbox{$\textstyle#1$}%
\vbox{\kern .33pt\hbox to\wd0{\kern 1pt\leaders\hrule height .33pt%
\hfill\kern 1pt}\kern 1pt\box0}}
{\setbox0=\hbox{$\scriptstyle#1$}%
\vbox{\kern .25pt\hbox to\wd0{\kern .8pt\leaders\hrule height .25pt%
\hfill\kern .8pt}\kern .8pt\box0}}
{\setbox0=\hbox{$\scriptscriptstyle#1$}%
\vbox{\kern .2pt\hbox to\wd0{\kern .6pt\leaders\hrule height .2pt%
\hfill\kern .6pt}\kern .6pt\box0}}}}

\def\det{\mathop{\rm det}\nolimits}
\def\End{\mathop{\rm End}\nolimits}

\def\Hom{\mathop{\rm Hom}\nolimits}
\def\id{\mathop{\rm id}\nolimits}

\def\Pic{\mathop{\rm Pic}\nolimits}

\def\SL{\mathop{\rm SL}\nolimits}

\def\Spec{\mathop{\rm Spec}\nolimits}
\def\Sym{\mathop{\rm Sym}\nolimits}

\def\tr{\mathop{\rm tr}\nolimits}

\def\isom{\buildrel\sim\over\longrightarrow}
\def\morphism#1{\buildrel#1\over\longrightarrow}
\def\isomorphism#1{\mathrel{\mathop{\longrightarrow}%
\limits^{#1}_{\raise0.5ex\hbox{$\scriptstyle\sim$}}}}
\def\leftisomorphism#1{\mathrel{\mathop{\longleftarrow}%
\limits^{#1}_{\raise0.5ex\hbox{$\scriptstyle\sim$}}}}

\def\injlim{\mathop{\vtop{\offinterlineskip\halign{##\cr
 \hfil\rm lim\hfil\cr\noalign{\kern.1ex}\rightarrowfill\cr
 \noalign{\kern-.4ex}\cr}}}}
\def\projlim{\mathop{\vtop{\offinterlineskip\halign{##\cr
 \hfil\rm lim\hfil\cr\noalign{\kern.1ex}\leftarrowfill\cr
 \noalign{\kern-.4ex}\cr}}}}

\def\blank{\mkern12mu}

\def\textfrac#1/#2{{\textstyle{#1\over#2}}}

\def\commdiag#1{{
\def\rightar{\longrightarrow}
\def\downar{\big\downarrow}
\def\upar{\big\uparrow}
\def\hookrightar{\lhook\joinrel\longrightarrow}
\def\rightlabel##1{\rlap{$\scriptstyle##1$}}
\def\leftlabel##1{\llap{$\scriptstyle##1$}}
\vcenter{\baselineskip=0pt \lineskiplimit=0pt \lineskip=6pt
\halign{\hfil$##$\hfil&&\enskip\hfil$##$\hfil\crcr
#1\crcr}}}}

\def\relativediag#1#2#3#4#5#6{\vcenter{\baselineskip=3ex \halign{
\hfil$##$&$##$&$##$\hfil\cr
#1\quad& \hfilneg\buildrel#2\over\longrightarrow\hfilneg& \quad#3\cr
\lower1ex\llap{$\scriptstyle#4\hskip-1ex$}\searrow& \quad&
\swarrow\lower1ex\rlap{$\hskip-1ex\scriptstyle#5$} \cr
& \hfilneg#6\hfilneg&\cr}}}

\def\trianglediag#1#2#3#4#5#6{\vcenter{\baselineskip=3ex \halign{
\hfil$##$&$##$&$##$\hfil\cr
#1\quad& \hfilneg\buildrel#2\over\longrightarrow\hfilneg& \quad#3\cr
\lower1ex\llap{$\scriptstyle#6\hskip-1ex$}\nwarrow& \quad&
\swarrow\lower1ex\rlap{$\hskip-1ex\scriptstyle#4$} \cr
& \hfilneg#5\hfilneg&\cr}}}

\def\correspondence#1#2#3#4#5{\vcenter{\baselineskip=3ex \halign{
\hfil$##$&\hfil$##$&$##$&$##$\hfil&$##$\hfil\cr
&&\hfilneg#1\hfilneg\cr
\smash{\raise1.5ex\hbox{$\scriptstyle#2$}}&\swarrow&&\searrow&
\smash{\raise1.5ex\hbox{$\scriptstyle#3$}}\cr
#4&&&&#5\cr}}}

\outer\def\abstract#1{{\narrower\narrower\noindent{\it Abstract\/.\enspace}#1\par}}

\outer\def\chapter#1. #2{\vfill\eject\ifodd\pageno\else\null\vfill\eject\fi
\sectioncount=0\edef\currentlabel{#1}%
\baselineskip=\chapterbaselineskip
\hangindent\chapterindent\hangafter1{\chapterfont\noindent
\hbox to\chapterindent{#1\hfill}#2}%
\bigskip\nobreak\normalbaselines\parindent=0pt\firstpartrue}

\newif\iffirstpar
\everypar{\iffirstpar\parindent=\normalparindent\firstparfalse\fi}

\outer\def\section#1{\advance\sectioncount by1
  \subcount=0 \subsectioncount=0 \eqcount=0
  \edef\currentlabel{\number\sectioncount}%
  \vskip 0pt plus.3\vsize\penalty-250
  \vskip 0pt plus-.3\vsize\bigskip\vskip\parskip
  \baselineskip=\sectionbaselineskip
  \hangindent\sectionindent\hangafter1
  {\sectionfont
  \noindent\hbox to\sectionindent{\number\sectioncount.\hfill}#1}%
  \medskip\nobreak\normalbaselines\parindent=0pt\firstpartrue}

\outer\def\unnumberedsection#1{
  \vskip 0pt plus.3\vsize\penalty-250
  \vskip 0pt plus-.3\vsize\bigskip\vskip\parskip
  \baselineskip=\sectionbaselineskip\noindent
  {\sectionfont #1}%
  \medskip\nobreak\normalbaselines\parindent=0pt\firstpartrue}

\outer\def\subsection#1{\advance\subsectioncount by1%
  \edef\currentlabel{\number\sectioncount.\number\subsectioncount}%
  \medbreak\hangindent\subsectionindent\hangafter1
  {\subsectionfont\noindent\hbox to\subsectionindent{%
  \number\sectioncount.\number\subsectioncount.\hfill}#1}%
  \smallskip\nobreak\parindent=0pt\firstpartrue}

\newwrite\auxfile

\newcount\sectioncount \sectioncount=0
\newcount\subsectioncount 
\newcount\subcount 
\newcount\eqcount 

\def\subno{\global\advance\subcount by1\relax
  \number\sectioncount.\number\subcount
  \xdef\currentlabel{\number\sectioncount.\number\subcount}}
\outer\def\proclaim #1. #2\par{\medbreak
  \noindent{\bf#1~\subno.\enspace}{\sl#2\par}%
  \ifdim\lastskip<\medskipamount \removelastskip\penalty55\medskip\fi}
\outer\def\proclaimx #1 (#2). #3\par{\medbreak
  \noindent{\bf#1~\subno\ \rm (#2).\enspace}{\sl#3\par}%
  \ifdim\lastskip<\medskipamount \removelastskip\penalty55\medskip\fi}

\newdimen\algindent
\def\plusindent{\advance\algindent by \parindent}
\def\minusindent{\advance\algindent by-\parindent}


\newcount\algstepcount

\long\def\algorithm (#1). #2\endalgorithm{\medbreak
  \algindent=0pt%
  \algstepcount=0%
  \noindent{\bf Algorithm~\subno} {\sl (#1)\/. \rm #2}\par
  \ifdim\lastskip<\medskipamount \removelastskip\penalty55\medskip\fi}%
\def\step{\advance\algstepcount by1
\edef\currentlabel{\number\algstepcount}%
\smallskip\hangindent\parindent
\advance\hangindent by\algindent\indent
\llap{{\bf \the\algstepcount.}\enspace}\kern\algindent}

\def\analysis{\noindent{\it Analysis\/}.\enspace}
\def\endanalysis{\nobreak\hfill\quad$\diamond$\par\medbreak}


\def\labeldef#1#2{\expandafter\gdef\csname L@#1\endcsname{#2}}
\def\label#1{%
  \expandafter\xdef\csname L@#1\endcsname{\currentlabel}%
  \write\auxfile{\string\labeldef{#1}{\csname L@#1\endcsname}}}
\def\ref#1{\expandafter\ifx\csname L@#1\endcsname\relax
  \message{Undefined label `#1'}??\else
  \csname L@#1\endcsname\fi}

\def\eqnumber#1{\global\advance\eqcount by1\relax
  \eqno(\number\sectioncount.\number\eqcount)%
  \expandafter\xdef\csname E@#1\endcsname{%
    \number\sectioncount.\number\eqcount}}
\def\eqref#1{\expandafter\ifx\csname E@#1\endcsname\relax
  \message{Undefined equation `#1'}??\else
  (\csname E@#1\endcsname)\fi}

\newcount\refcount \refcount=0
\def\citedef#1#2{\expandafter\gdef\csname C@#1\endcsname{#2}}
\def\citepage#1{%
  \expandafter\ifx\csname Rlastcited@#1\endcsname\relax
    \expandafter\xdef\csname Rcited@#1\endcsname{\the\pageno}%
  \else
     \edef\pn{\the\pageno}%
     \expandafter\ifx\csname Rlastcited@#1\endcsname\pn
     \else
     \expandafter\xdef\csname Rcited@#1\endcsname{%
        \csname Rcited@#1\endcsname, \the\pageno}%
     \fi
  \fi
  \expandafter\xdef\csname Rlastcited@#1\endcsname{\the\pageno}}
\def\cite#1{\expandafter\ifx\csname C@#1\endcsname\relax
  \message{Undefined reference `#1'}\citedef{#1}{??}\fi
  \citepage{#1}%
  [\csname C@#1\endcsname]}
\def\citex#1#2{\expandafter\ifx\csname C@#1\endcsname\relax
  \message{Undefined reference `#1'}\citedef{#1}{??}\fi
  \citepage{#1}%
  [\csname C@#1\endcsname, #2]}
\def\reference#1#2\par{\advance\refcount by 1%
  \expandafter\edef\csname C@#1\endcsname{\the\refcount}%
  \write\auxfile{\string\citedef{#1}{\csname C@#1\endcsname}}%
  \item{[\csname C@#1\endcsname]}#2%
  \expandafter\ifx\csname Rcited@#1\endcsname\relax
  \message{Warning: reference `#1' not used}%
  \else{\it \csname Rcited@#1\endcsname}\fi}

\newwrite\indfile
\newif\ifmakeind \makeindfalse

\def\writeind#1#2{\write\indfile{\string\ind{#1}{#2}}}
\def\index#1{\ifmakeind\writeind{#1}{\the\pageno}\fi}
\def\writeindsub#1#2#3{\write\indfile{\string\ind{#1}\string\sub{#2}{#3}}}
\def\indexsub#1#2{\ifmakeind\writeindsub{#1}{#2}{\the\pageno}\fi}

\newif\ifauxexists
\immediate\openin0=\jobname.aux
\ifeof 0
  \auxexistsfalse
\else
  \auxexiststrue
\fi
\immediate\closein0
\ifauxexists
  \input \jobname.aux
\else
  \message{No file `\jobname.aux'}
\fi
\openout\auxfile=\jobname.aux

\tenpoint
\normalbaselines
\rm

\def\Alb{\mathop{\rm Alb}\nolimits}

\def\div{\mathop{\rm div}\nolimits}
\def\Div{\mathop{\rm Div}\nolimits}
\def\Eff{\mathop{\rm Ef\/f}\nolimits}
\def\F{{\bf F}}

\def\Frob{{\rm F}}
\def\G{{\bf G}}
\def\genus_#1{g_{#1}^{\null}}  
\def\glob(#1,#2){\Gamma({#1},{#2})}  
\def\gl(#1,#2){\Gamma({#2})}         
\def\Gr{\mathop{\rm Gr}\nolimits}

\def\Irr{\mathop{\rm Irr}\nolimits}
\def\L{{\cal L}}
\def\Norm{{\rm N}}
\def\O{{\cal O}}
\def\P{{\bf P}}
\def\PrimeDiv{\mathop{\rm PDiv}\nolimits}

\def\Q{{\bf Q}}

\def\shHom{\mathop{\cal H\mit om}\nolimits}
\def\supp{\mathop{\rm supp}\nolimits}

\def\Ver{\mathop{\rm Ver}\nolimits}
\def\X{{\rm X}}
\def\Z{{\bf Z}}
\def\Zeta{{\rm Z}}

\def\PhD{Ph.\thinspace D.\ }

\centerline{\titlefont Computing in Picard groups}
\medskip
\centerline{\titlefont of projective curves over finite fields}
\bigskip
\centerline{Peter Bruin}
\par

\bigskip

\abstract{We give algorithms for computing with divisors on projective
curves over finite fields, and with their Jacobians, using the
algorithmic representation of projective curves developed by
Khuri-Makdisi.  We show that various desirable operations can be
performed efficiently in this setting: decomposing divisors into prime
divisors; computing pull-backs and push-forwards of divisors under
finite morphisms, and hence Picard and Albanese maps on Jacobians;
generating uniformly random divisors and points on Jacobians;
computing Frobenius maps and Kummer maps; and finding a basis for the
$l$-torsion of the Picard group, where $l$ is a prime number different
from the characteristic of the base field.}

\vskip.4cm

\unnumberedsection{Introduction}

In \cite{Khuri-Makdisi: Linear algebra algorithms}
and~\cite{Khuri-Makdisi}, K. Khuri-Makdisi developed efficient
algorithms for computing with divisors on projective curves over
arbitrary fields.  The goal of this article is to show that for curves
over finite fields in Khuri-Makdisi's algorithmic representation, one
can compute Frobenius morphisms and Frey--R\"uck pairings, pick
uniformly random rational points on curves and their Jacobians (given
the zeta function of the curve), perform various other operations
specific to curves over finite fields, and compute Picard and Albanese
maps induced by certain finite morphisms between curves.

The curves we consider are complete, smooth and geometrically
connected curves over a field~$k$.  For now we assume $k$ is an
arbitrary field; later we assume it to be finite.  The basic idea is
to describe such a curve~$X$ using a projective embedding via a very
ample line bundle~$\L$.  The curve is then represented by means of the
finite $k$-algebra obtained by taking the quotient of the homogeneous
coordinate ring of~$X$ by the ideal generated by homogeneous elements
of sufficiently large degree.  Divisors on~$X$ are represented as
subspaces of the $k$-vector spaces of global sections of suitable
powers of the line bundle~${\cal L}$.  Using this representation of
the curve and of divisors on it, Khuri-Makdisi~\cite{Khuri-Makdisi:
Linear algebra algorithms} has given algorithms for computing with
divisors and elements of the Picard group.  Taking advantage of some
improvements to this basic idea, described in~\cite{Khuri-Makdisi},
his algorithms are currently the fastest known algorithms for general
curves, asymptotically as the genus increases and measured in
operations in~$k$.

The algorithms presented in this paper are relevant for computations
with curves of large genus over finite fields.  The author's interest
in these was raised by algorithms for explicitly computing
coefficients of modular forms.  In \cite{compcoefs} and the
forthcoming book~\cite{book}, Couveignes, Edixhoven and others
describe an algorithm for computing coefficients of modular forms for
the group~$\SL_2(\Z)$.  In the author's forthcoming
thesis~\cite{thesis}, their methods are generalised to modular forms
for groups of the form~$\Gamma_1(n)$.  The method used in each case is
to compute two-dimensional modular Galois representations over finite
fields.  The basic problem is to find explicit realisations of group
schemes over~$\Q$ of the form $J[{\frak m}]$, where $J$ is the
Jacobian of a modular curve and ${\frak m}$ is a maximal ideal of the
corresponding Hecke algebra.  The approach taken is to approximate the
coefficients of certain polynomials defining such group schemes,
either over the complex numbers or modulo sufficiently many small
prime numbers.  The complex method has already been used by
Bosman~\cite{Bosman} for actual computations.  The alternative method
using finite fields was described by Couveignes in~\cite{compcoefs}
for the modular curves~$\X_1(5l)$, where $l$ is a prime number.  The
computations in this case can be done using (singular) plane models
for these curves.  For a more general modular curve~$X$, it seems
natural to take an embedding of~$X$ as a smooth curve in a
higher-dimensional projective space, using the line bundle of modular
forms of weight~2.  Using the technique of modular
symbols~\cite{Stein}, one can compute $q$-expansions of these modular
forms, as well as the zeta function of~$X$.  This immediately gives a
representation of~$X$ that can be used for Khuri-Makdisi's algorithms,
without having to write down equations.

\goodbreak
The paper is organised as follows.  In the preliminary
Section~\ref{sec:finite-algebras} we consider some computational
problems related to finite algebras over a field; these are needed in
the other two sections.  In Section~\ref{sec:computing-with-divisors}
we recall Khuri-Makdisi's algorithms for projective curves over
arbitrary base fields, and we describe a number of extensions.  Some
of our algorithms require that we are able to efficiently compute
primary decompositions of finite $k$-algebras.  This condition is
fulfilled, for example, if $k$ is a finite field or a number field.
We give algorithms for the following computational problems:
\smallskip
\item{(1)} finding the decomposition of a divisor as a
linear combination of prime divisors;
\smallskip
\item{(2)} computing pull-backs and push-forwards of divisors
under finite morphisms;
\smallskip
\item{(3)} computing Picard and Albanese maps induced by finite
morphisms of curves.
\smallskip
\noindent We also consider some more technical problems that are
needed for the rest of the paper.  In
Section~\ref{sec:curves-over-finite-fields} we describe the rest of
our algorithms, which are specific to curves over finite fields.
These are the following:
\smallskip
\item{(1)} computing the Frobenius map on points of the curve, and of
its Jacobian, that are defined over finite extensions of the base
field;
\smallskip
\item{(2)} generating uniformly random effective divisors of a given
degree, and uniformly random points of the Jacobian, if the zeta
function of the curve is known;
\smallskip
\item{(3)} computing Frey--R\"uck pairings on the Jacobian.
\smallskip\noindent
By combining the above methods, we also show that the methods of
Couveignes~\cite{Couveignes: Linearizing torsion classes} for
computing Kummer maps of order~$l$ and for finding a basis for the
$l$-torsion of the Picard group, where $l$ is a prime number different
from the characteristic of the base field, can be extended to our
situation, again under the assumption that we know the zeta function
of the curve.
\medskip

\remarks (1)\enspace When the field~$k$ is finite, measuring the
running time in field operations is essentially the same as measuring
it in bit operations.  However, if $k$ is a number field, it is
impossible to avoid numerical explosion of the data describing the
divisors during computations, so that the running time in bit
operations is much worse than that counted in bit operations.  Using
lattice reduction algorithms to reduce the size of the data between
operations should not be expected to solve this problem; see
Khuri-Makdisi~\citex{Khuri-Makdisi}{page~2214}.
\smallskip\noindent
(2)\enspace Many of the algorithms we describe are probabilistic.  All
of these are of the {\it Las Vegas\/} type, meaning that the running
time depends on certain random data generated during the execution of
the algorithm, but that the outcome is guaranteed to be correct.  The
epithet {\it Las Vegas\/} distinguishes such algorithms from those of
the {\it Monte Carlo\/} type, where the randomness influences the
correctness of the outcome instead of the running time.
\smallskip\noindent
(3)\enspace The algorithms mentioned in this paper have a running time
that is bounded by some polynomial in various quantities that are
indicated in each case.  Obtaining more detailed estimates should not
be difficult, but has at the time of writing not yet been done.

\acknowledgements I would like to thank Johan Bosman, Claus Diem,
Bas Edixhoven, Kamal Khuri-Makdisi and Hendrik Lenstra for useful
conversations and correspondence on topics related to this paper.

\section{Algorithms for computing with finite algebras}

\label{sec:finite-algebras}

In this section we describe techniques for solving two computational
problems about finite algebras over a field.  The first is how to find
the primary decomposition of such an algebra; the second is how to
reconstruct such an algebra from a certain kind of bilinear map beteen
modules over it.

The algebras to which we are going to apply these techniques in the
next section are of the form~$\glob(E,\O_E)$, where $E$ is an
effective divisor on a smooth curve over~$k$.  In this section,
however, we place ourselves in the more general setting of arbitrary
finite commutative $k$-algebras.

\subsection{Primary decomposition and radicals}

\label{primary-decomposition-radicals}

Let $k$ be a field with the following two properties:
\smallskip
\item{(1)} $k$ is perfect;
\smallskip
\item{(2)} we have a (probabilistic) algorithm to factor polynomials
$f\in k[x]$ that takes an (expected) number of operations in~$k$ that
is bounded by a polynomial in the degree of~$f$.
\smallskip
\noindent
For such a field~$k$ there exist (probabilistic) algorithms to find
the primary decomposition of a finite commutative $k$-algebra~$A$ that
finish in an (expected) number of operations in~$k$ that is bounded by
a polynomial in~$[A:k]$.  Such algorithms have been known for some
time, but do not seem to be easily available in published form; see
Khuri-Makdisi's preprint~\citex{Khuri-Makdisi}{draft version~2, \S7}.
For an algorithm to find the primary decomposition of arbitrary (not
necessarily commutative) finite algebras over {\it finite\/} fields,
see Eberly and Giesbrecht~\cite{Eberly-Giesbrecht}.

\subsection{Reconstructing an algebra from a perfect bilinear map}

\label{reconstructing-algebra-from-bilinear-map}

Let $A$ be a commutative ring.  If $M$, $N$ and~$O$ are free
$A$-modules of rank~one and
$$
\mu\colon M\times N\to O
$$
is an $A$-bilinear map, we say that $\mu$ is {\it perfect\/} if it
induces an isomorphism
$$
M\otimes_A N\isom O
$$
of free $A$-modules of rank~1.  

Now let $k$ be a field, and let a finite commutative $k$-algebra $A$
be specified implicitly in the following way.  We are given $k$-vector
spaces $M$ and~$N$ of the same finite dimension, together with a
$k$-bilinear map
$$
\mu\colon M\times N\to O
$$
We assume there exists a commutative $k$-algebra~$A$ such that $M$,
$N$ and~$O$ are free $A$-modules of rank~1 and $\mu$ is a perfect
$A$-bilinear map.  The following observation implies that $A$ is the
{\it unique\/} $k$-algebra with this property, and also shows how to
compute $A$ as a subalgebra of~$\End_k M$, provided we are able to
find a generator of~$N$ as an $A$-module.  As could be expected, the
roles of $M$ and~$N$ can be interchanged.

\proclaim Lemma. In the above situation, let $g$ be a generator of the
$A$-module $N$.  The ring homomorphism $A\to\End_k M$ sending $a$ to
multiplication by~$a$ is, as an $A$-linear map, the composition of
$$
\eqalignno{
A&\isom N\cr
a&\longmapsto ag\cr
\noalign{\noindent and}
N&\longrightarrow\End_k M\cr
n&\longmapsto\mu(\blank,g)^{-1}\circ\mu(\blank,n).}
$$
In particular, the image of~$A$ in~$\End_kM$ equals the image of the
second map.

\label{endomorphism-algebra-lemma}

\proof This is a straightforward verification.\endproof

In the case where $k$ is a finite field, a way to find a generator
for~$N$ as an $A$-module is simply to pick random elements $g\in N$
until we find one that generates $N$.  Since $\mu$ is perfect,
checking whether $g$ generates $N$ comes down to checking whether
$\mu(\blank,g)\colon M\to O$ is an isomorphism.  In particular, we can
do this without knowing $A$.

To get a reasonable expected running time for this approach, we need
to ensure that $N$ contains sufficiently many elements that generate
it as an $A$-module.  Since $N$ is free of rank~1, the number of
generators equals the number of units in~$A$.  Let us therefore
estimate under what conditions a random element of~$A$ is a unit with
probability at least~1/2.  Write $d$ for the degree of $A$ over~$k$.
Decomposing $A$ into a product of finite local $k$-algebras, and
noting that the proportion of units in a finite local $k$-algebra is
equal to the proportion of units in its residue class field, we see
that
$$
{\#A^\times\over\#A}\ge{(\#k^\times)^d\over\#k^d}
=\left(1-{1\over\#k}\right)^d;
$$
equality occurs if and only if $A$ is a product of $d$ copies of~$k$.
Now it is not hard to show that
$$
\#k\ge 2d\;\Longrightarrow\;
\left(1-{1\over\#k}\right)^d\ge{1\over2}.
$$
Taking a finite extension~$k'$ of~$k$ of cardinality at least $2d$, we
therefore see that a random element of~$A_{k'}$ is a unit with
probability at least $1/2$.  There are well-known algorithms to
generate such an extension, such as that of Rabin~\cite{Rabin}, which
runs in probabilistic polynomial time and simply tries random
polynomials until it finds one that is irreducible, and the
deterministic algorithm of Adleman and Lenstra~\cite{Adleman-Lenstra},
which is only known to run in polynomial time under the generalised
Riemann hypothesis.

\algorithm (Reconstruct an algebra from a bilinear map).  Let $k$ be a
finite field, let $A$ be a finite $k$-algebra, and let
$$
\mu\colon M\times N\to O
$$
be a perfect $A$-bilinear map between free $A$-modules of rank~1.
Given the coefficients of~$\mu$ with respect to some $k$-bases of $M$,
$N$ and~$O$, this algorithm outputs a $k$-basis for the image of~$A$
in~$\End_k M$, consisting of matrices with respect to the given basis
of~$M$.

\step Choose an extension $k'$ of~$k$ of
degree~$\left\lceil{\log\max\{2[A:k],q\}\over\log q}\right\rceil$.
Let $M'$, $N'$, $O'$ and~$\mu'$ denote the base extensions of~$M$,
$N$, $O$ and~$\mu$ to~$k'$.

\step Choose a uniformly random element $g\in N'$.

\step Check whether $\mu'(\blank,g)\colon M'\to O'$ is an isomorphism;
if not, go to step~2.

\step For $n$ ranging over a $k'$-basis of $N'$, compute the
endomorphism
$$
a_n=\mu'(\blank,g)^{-1}\circ\mu'(\blank,n)\in\End_{k'}M'.
$$
Let $A'\subseteq\End_{k'}M'$ denote the $k'$-span of the $a_n$.

\step Output a basis for the $k$-vector space $\End_k M\cap A'$.

\endalgorithm

\analysis It follows from Lemma~\ref{endomorphism-algebra-lemma} that
$A'$ equals the image of~$k'\otimes_k A$ in~$\End_{k'} M$.  This
implies that the basis returned by the algorithm is indeed a $k$-basis
for the image of~$A$ in~$\End_k M$.  Because of the choice of~$k'$,
steps 2 and~3 are executed at most twice on average.  It is therefore
clear that the expected running time of the algorithm is polynomial in
$[A:k]$ and $\log\#k$.\endanalysis

\label{bilinear-map-to-algebra-finite}

If $k$ is infinite (or finite and sufficiently large), we have the
following variant.  Let $\Sigma$ be a finite subset of~$k$, and let
$V$ be a $k$-vector space of dimension~$d$ with a given basis
$v_1$, \dots, $v_d$.  Consider the set
$$
V_\Sigma=\{\sum_{i=1}^d\sigma_i v_i\mid
  \sigma_1,\ldots,\sigma_d\in\Sigma\}
$$
of {\it $\Sigma$-linear combinations\/} of $v_1$, \dots, $v_n$.
Choosing the $\sigma_i$ uniformly randomly in~$\Sigma$, we get the
uniform distribution on~$V_\Sigma$.  If $H_1$, \dots, $H_l$ are proper
linear subspaces of~$V$, then a uniformly random element of~$V_\Sigma$
lies in at least one of the $H_i$ with probability at most
$l/\#\Sigma$.  Now if $A$ is a finite commutative $k$-algebra, it
contains at most $[A:k]$ maximal ideals.  This implies that if
$\Sigma$ is a finite subset of~$k$ with $\#\Sigma\ge 2[A:k]$, then a
$\Sigma$-linear combination of any $k$-basis of~$A$ is a unit with
probability at least $1/2$.  This leads to the following variant of
Algorithm~\ref{bilinear-map-to-algebra-finite}.

\algorithm (Reconstruct an algebra from a bilinear map).  Let $k$ be a
field, let $A$ be a finite $k$-algebra, and let
$$
\mu\colon M\times N\to O
$$
be a perfect $A$-bilinear map between free $A$-modules of rank~1.
Suppose that we can pick uniformly random elements of some subset
$\Sigma$ of~$k$ with $\#\Sigma\ge 2[A:k]$.  Given the coefficients
of~$\mu$ with respect to some $k$-bases of $M$, $N$ and~$O$, this
algorithm outputs a $k$-basis for the image of~$A$ in~$\End_k M$,
consisting of matrices with respect to the given basis of~$M$.

\step Choose a uniformly random $\Sigma$-linear combination $g$ of
the given basis of~$N$.

\step Check whether $\mu(\blank,g)\colon M\to O$ is an isomorphism;
if not, go to step~2.

\step For $n$ ranging over a $k$-basis of $N$, compute the
endomorphism
$$
a_n=\mu(\blank,g)^{-1}\circ\mu(\blank,n)\in\End_k M,
$$
and output the $a_n$.

\endalgorithm

\label{bilinear-map-to-algebra-infinite}

\analysis This works for the same reason as
Algorithm~\ref{bilinear-map-to-algebra-finite}.\endanalysis

Let us sketch how to solve the problem if $k$ is an arbitrary field.
Let $p$ be the characteristic of~$k$.  If $p=0$ or $p\ge 2[A:d]$, we
can apply Algorithm~\ref{bilinear-map-to-algebra-infinite} with
$\Sigma=\{0,1,\ldots,2[A:d]-1\}$.  Otherwise, we consider the
subfield~$k_0$ of~$k$ generated by the coefficients of the
multiplication table of~$A$ over~$k$.  Then $A$ is obtained by base
extension to~$k$ of the finite $k_0$-algebra $A_0$ defined by the same
multiplication table.  We can check whether $k_0$ is a finite field
with $\#k_0<2d$ by checking whether each coefficient of the
multiplication table satisfies a polynomial of small degree.  If this
is the case, then we compute an $\F_p$-basis and multiplication table
for~$k_0$ and apply Algorithm~\ref{bilinear-map-to-algebra-finite} to
$A_0$ over~$k_0$.  Otherwise we obtain at some point a finite
subset~$\Sigma$ of~$k$, with $\#\Sigma\ge 2d$, consisting of
polynomials in the coefficients of the multiplication table.  We then
apply Algorithm~\ref{bilinear-map-to-algebra-infinite} to $A$ over~$k$
with this $\Sigma$.

\section{Computing with divisors on a curve}

\label{sec:computing-with-divisors}

In this section we describe a collection of algorithms, developed by
Khuri-Makdisi in \cite{Khuri-Makdisi: Linear algebra algorithms}
and~\cite{Khuri-Makdisi}, that allow us to compute efficiently with
divisors on a curve over a field.  In particular, we will describe
algorithms for computing in the Picard group of a curve.  Many of the
results of this section can be found in \cite{Khuri-Makdisi: Linear
algebra algorithms} and~\cite{Khuri-Makdisi};
however, \SS\ref{functoriality-for-finite-morphisms},
\ref{normalised-representatives} and~\ref{picard-albanese} seem to be
new.

\subsection{Representing the curve}

\label{representation-curve}

Let $X$ be a complete, smooth, geometrically connected curve over a
field~$k$.  We fix a line bundle~$\L$ on~$X$ such that
$$
\deg\L\ge 2g+1.
$$
Then $\L$ is very ample (see for example
Hartshorne~\citex{Hartshorne}{IV, Corollary~3.2(b)}), so it gives rise
to a closed immersion
$$
i_\L\colon X\to\P\glob(X,\L)
$$
into a projective space of dimension $\deg\L-g$.  (We write $\P V$ for
the projective space of hyperplanes in a $k$-vector space~$V$.)  The
assumption that $\deg\L\ge 2g+1$ implies moreover that the
multiplication maps
$$
\mu_{i,j}\colon\glob(X,\L^{\otimes i})\otimes_k
\glob(X,\L^{\otimes j})\longrightarrow
\glob(X,\L^{\otimes(i+j)}).
$$
are surjective for all $i,j\ge0$, or equivalently that the
embedding~$i_\L$ is projectively normal.  This is a classical theorem
of Castelnuovo, Mattuck and Mumford; see for example
Lazarsfeld~\citex{Lazarsfeld}{\S1.1}.


We write $S_X$ for the homogeneous coordinate ring of~$X$ with respect
to the embedding~$i_\L$.  By the fact that $i_\L$ is
projectively normal, we have a canonical isomorphism
$$
S_X\isom\bigoplus_{i\ge0}\glob(X,\L^{\otimes i})
$$
of graded $k$-algebras; see Hartshorne~\citex{Hartshorne}{Chapter~II,
Exercise~5.14}.  It turns out that to be able to compute with divisors
on~$X$ we do not need to know the complete structure of this graded
algebra.  For all $h\ge0$ we define the finite graded $k$-algebra
$S_X^{(h)}$ as $S_X$ modulo the ideal generated by homogeneous
elements of degree greater than~$h$.  The above isomorphism shows that
specifying $S_X^{(h)}$ is equivalent to giving the $k$-vector spaces
$\glob(X,\L^{\otimes i})$ for $1\le i\le h$ together with the
multiplication maps~$\mu_{i,j}$ between them for $i+j\le h$.

When we speak of a {\it projective curve\/}~$X$ in the remainder of
this section, we will assume without further mention that $X$ is a
complete, smooth and geometrically connected curve of genus~$g\ge0$,
and that a line bundle~$\L$ of degree at least~$2g+1$ has been chosen.
We will often write $\L_X$ for this line bundle and $\genus_X$ for the
genus of~$X$ to emphasise that they are part of the data.

In the algorithms in this section, the curve~$X$ is part of the input
in the guise of the graded $k$-algebra $S_X^{(h)}$ for some
sufficiently large $h$.  A lower bound for~$h$ is specified in each
case.  One way to specify the multiplication in~$S_X^{(h)}$ is to fix
a basis for each of the spaces $\glob(X,\L^{\otimes i})$, and to give
the matrices for multiplication with each basis element.  However, as
Khuri-Makdisi explains in~\cite{Khuri-Makdisi}, a more efficient
representation is to choose a trivialisation of~$\L$ (and hence of its
powers) over an effective divisor of sufficiently large degree or,
even better, at sufficiently many distinct rational points of~$X$, so
that the multiplication maps can be computed pointwise.

\remarks (1)\enspace The integers $g$ and~$\deg\L$ can of course be
stored as part of the data describing~$X$.  However, they can also be
extracted from the dimensions of the $k$-vector spaces $\glob(X,\L)$
and~$\glob(X,\L^{\otimes2})$, and hence from~$S_X^{(2)}$; this follows
easily from the Riemann--Roch formula.
\smallskip\noindent
(2)\enspace If the degree of~$\L$ is at least $2g+2$, then the
homogeneous ideal defining the embedding~$i_\L$ is generated by
homogeneous elements of degree~2, according to a theorem of Fujita and
Saint-Donat; see Lazarsfeld~\citex{Lazarsfeld}{\S1.1}.  This makes it
possible to deduce equations for~$X$ from the $k$-algebra $S_X^{(2)}$.
However, we will not need to do this.
\smallskip\noindent
(3)\enspace The representation of curves described by Khuri-Makdisi
in \cite{Khuri-Makdisi: Linear algebra algorithms}
and~\cite{Khuri-Makdisi} is especially suited for modular curves.
Namely, we can represent a modular curve~$X$ using the projective
embedding given by a line bundle of modular forms, and computing the
$k$-algebra $S_X^{(h)}$ for a given $h$ comes down to computing
$q$-expansions of modular forms of a suitable weight to sufficient
precision.  This can be done using modular symbols; see
Stein~\cite{Stein}.  If the modular curve has at least 3 cusps (which
is the case, for example, for ${\rm X}_1(n)$ for all $n\ge5$), then we
can restrict ourselves to modular forms of weight~2, for which the
formalism of modular symbols is particularly
simple~\citex{Stein}{Chapter~3}.

\subsection{Representing divisors}

\label{representation-divisors}

Let $X$ be a projective curve of genus~$g$ in the sense
of~\S\ref{representation-curve}, and let $\L$ be the line bundle
giving the projective embedding of~$X$.  To represent divisors on~$X$,
it is enough to consider effective divisors, since an arbitrary
divisor can be represented by a formal difference of two effective
divisors.

Let $D$ be an effective divisor on~$X$ such that $\L(-D)$ is generated
by global sections.  In terms of the projective embedding, this means
that $D$ is the intersection of $X$ and a linear subvariety
of~$\P\glob(X,\L)$, or equivalently that $D$ is defined by a system of
linear equations.  Such a divisor can be represented as the subspace
$\glob(X,\L(-D))$ of~$\glob(X,\L)$ consisting of sections vanishing
on~$D$.  The codimension of $\glob(X,\L(-D))$ in~$\glob(X,\L)$ is
equal to the degree of~$D$.

A sufficient condition for the line bundle~$\L(-D)$ to be
generated by global sections is
$$
\deg D\le\deg\L-2g;
\eqnumber{degree-bound}
$$
see for example Hartshorne~\citex{Hartshorne}{IV, Corollary~3.2(a)}.
However, we note that in general not every subspace of codimension at
most~$\deg\L-2g$ is of the form $\glob(X,\L(-D))$ for an
effective divisor~$D$ of the same degree.

\remark This way of representing divisors comes down (at least for
divisors of degree $d\le\deg\L-2g$) to embedding the $d$-th
symmetric power of~$X$ into the Grassmannian variety parametrising
subspaces of codimension~$d$ in~$\glob(X,\L)$ and viewing
divisors of degree~$d$ as points on this Grassmannian variety.

It will often be necessary to consider divisors~$D$ of degree larger
than the bound $\deg\L-2g$ of~\eqref{degree-bound}.  In such cases we
can represent $D$ as a subspace of $\glob(X,\L^{\otimes i})$ for $i$
sufficiently large such that
$$
\deg D\le i\deg\L-2g,
\eqnumber{general-degree-bound}
$$
provided of course that we know $S_X^{(h)}$ for some $h\ge i$.

Khuri-Makdisi's algorithms rest on the following two results.  The
first is a generalisation of the theorem of Castelnuovo, Mattuck and
Mumford mentioned above.  It says in effect that to compute the
space of global sections of the tensor product of two line bundles of
sufficiently large degree, it is enough to multiply global sections of
those line bundles.

\proclaimx Lemma (Khuri-Makdisi~\citex{Khuri-Makdisi: Linear algebra
algorithms}{Lemma 2.2}). Let $X$ be a complete, smooth, geometrically
connected curve of genus~$g$ over a field~$k$, and let ${\cal M}$
and~${\cal N}$ be line bundles on~$X$ whose degrees are at least
$2g+1$.  Then the canonical $k$-linear map
$$
\glob(X,{\cal M})\otimes_k\glob(X,{\cal N})\longrightarrow
\glob(X,{\cal M}\otimes_{\O_X}{\cal N})
$$
is surjective.

\label{multiplication}

The second result shows how to find the space of global sections of a
line bundle that vanish on a given effective divisor, where this
divisor is represented as a subspace of global sections of a second
line bundle.

\proclaimx Lemma (Khuri-Makdisi~\citex{Khuri-Makdisi: Linear algebra
algorithms}{Lemma 2.3}). Let $X$ be a complete, smooth, geometrically
connected curve of genus~$g$ over a field~$k$, let ${\cal M}$
and~${\cal N}$ be line bundles on~$X$ such that ${\cal N}$ is
generated by global sections, and let $D$ be any effective divisor
on~$X$.  Then the inclusion
$$
\glob(X,{\cal M}(-D))\subseteq\bigl\{s\in\glob(X,{\cal M})\bigm|
s\glob(X,{\cal N})\subseteq\glob(X,{\cal M}\otimes{\cal N}(-D))
\bigr\}
$$
is an equality.

\label{division}

Thanks to these two lemmata, one can give algorithms to do basic
operations on divisors; see Khuri-Makdisi~\citex{Khuri-Makdisi: Linear
algebra algorithms}{\S3}.  For example, we can add, subtract and
intersect divisors of sufficiently small degree, and we can test
whether a given subspace of~$\glob(X,\L^{\otimes i})$ is of the form
$\glob(X,\L^{\otimes i}(-D))$ for some effective divisor~$D$.  See
also Algorithm~\ref{addflip} below for an example where
Lemmata \ref{multiplication} and~\ref{division} are used.

\subsection{Deflation and inflation}

\label{deflation-inflation}

An ingredient that Khuri-Makdisi uses in~\cite{Khuri-Makdisi} to speed
up the algorithms is {\it deflation\/} of subspaces.  Suppose we want
to compute the space $\glob(X,{\cal M}(-D))$ using
Lemma~\ref{division} in the case where ${\cal M}=\L^{\otimes i}$ and
${\cal N}=\L^{\otimes j}(-E)$ with $i$ and~$j$ positive integers and
where $D$ and~$E$ are effective divisors
satisfying \eqref{general-degree-bound}.  On the right-hand side of
the equality given by Lemma~\ref{division}, we may replace
$\glob(X,{\cal N})$ by any basepoint-free subspace; this is clear from
the proof of~\citex{Khuri-Makdisi: Linear algebra
algorithms}{Lemma~2.3}.  It turns out that there always exists such a
subspace of dimension $O(\log(\deg{\cal N}))$, and a subspace of
dimension~2 exists if the base field is either infinite or finite of
sufficiently large cardinality.  Moreover, one can efficiently find
such a subspace by random trial; see
Khuri-Makdisi~\citex{Khuri-Makdisi}{Proposition/Algorithm~3.7}.

Suppose we are given a basepoint-free subspace $W$
of~$\glob(X,\L^{\otimes i}(-D))$ for some $i$ and~$D$ such that
$\glob(X,\L^{\otimes i}(-D))$ is basepoint-free.  Then we can
reconstruct the complete space~$\glob(X,\L^{\otimes i}(-D))$ from~$W$.
This procedure is called {\it inflation\/}.  To describe how this can
be done, we first state the following slight generalisation of a
result of Khuri-Makdisi~\citex{Khuri-Makdisi}{Theorem 3.5(2)}.

\proclaim Lemma. Let $X$ be a complete, smooth, geometrically curve of
genus~$g$ over a field~$k$, and let ${\cal M}$ and~${\cal N}$ be line
bundles on~$X$.  Let $V$ be a non-zero subspace of~$\glob(X,{\cal
M})$, and let $D$ be the common divisor of the elements of~$V$.  If
the inequality
$$
-\deg{\cal M}+\deg{\cal N}+\deg D\ge 2g-1
$$
is satisfied, the canonical $k$-linear map
$$
V\otimes_k\glob(X,{\cal N})\longrightarrow
\glob(X,{\cal M}\otimes_{\O_X}{\cal N}(-D))
\eqnumber{multiplication-bpf}
$$
is surjective.

\label{lemma-mult-bpf}

\proof We note that ${\cal M}(-D)$ is basepoint-free by definition,
since we can view $V$ as a subspace of~$\glob(X,{\cal M}(-D))$ and the
elements of~$V$ have common divisor~0 as sections of~${\cal M}(-D)$.
We also note that $\deg{\cal M}\ge\deg D$.  Therefore the assumption
on the degrees of ${\cal M}$, ${\cal N}$ and~$D$ implies the
inequalities
$$
\deg{\cal N}\ge2g-1
$$
and
$$
\deg({\cal M}\otimes{\cal N}(-D))\ge 2g-1.
$$
After extending the field~$k$, we may assume it is infinite.  Then
there exist elements $s,t\in V$ with common divisor~$D$; see
Khuri-Makdisi~\citex{Khuri-Makdisi}{Lemma~4.1}.  The space
$$
s\glob(X,{\cal N})+t\glob(X,{\cal N})
$$
lies in the image of~\eqref{multiplication-bpf}, so it suffices to
show that
$$
\dim_k(s\glob(X,{\cal N})+t\glob(X,{\cal N}))=
\dim_k\glob(X,{\cal M}\otimes{\cal N}(-D)).
$$
Write
$$
\div s=D+E\quad\hbox{and}\quad\div s=D+F
$$
where $E$ and~$F$ are disjoint effective divisors.  Then we have
$$
\eqalign{
\dim_k(s\glob(X,{\cal N})+t\glob(X,{\cal N}))&=
2\dim_k\glob(X,{\cal N})-\dim_k(s\glob(X,{\cal N})
\cap t\glob(X,{\cal N}))\cr
&=2\dim_k\glob(X,{\cal N})\cr
&\qquad-\dim_k\glob(X,{\cal M}\otimes{\cal N}(-D-E-F))\cr
&=2\dim_k\glob(X,{\cal N})
-\dim_k\glob(X,{\cal M}^\vee\otimes{\cal N}(D)).}
$$
The last equality follows from the fact that multiplication by~$st$
induces an isomorphism
$$
{\cal M}^\vee(D)\isom{\cal M}(-D-E-F).
$$
Using the fact that the various line bundles have degrees at least
$2g-1$, we see that
$$
\eqalign{
\dim_k(s\glob(X,{\cal N})+t\glob(X,{\cal N}))&=
2(1-g+\deg{\cal N})-(1-g+\deg{\cal M}^\vee\otimes{\cal N}(D))\cr
&=1-g+\deg{\cal M}+\deg{\cal N}-\deg D\cr
&=\dim_k\glob(X,{\cal M}\otimes{\cal N}(-D)).}
$$
This finishes the proof.\endproof

We can now describe how to inflate a basepoint-free subspace $W$
of~$\glob(X,\L^{\otimes i}(-D))$.  Namely, we choose a positive
integer~$j$ such that
$$
(j-i)\deg\L+\deg D\ge 2g-1.
$$
By Lemma~\ref{lemma-mult-bpf} we can then compute
$\glob(X,\L^{\otimes(i+j)}(-D)) $ as the image of the bilinear map
$$
W\otimes_k\glob(X,\L^{\otimes j})\longrightarrow
\glob(X,\L^{\otimes(i+j)}).
$$
Then we compute
$$
\glob(X,\L^{\otimes i}(-D))=\bigl\{s\in\glob(X,\L^{\otimes i})\bigm|
s\glob(X,\L^{\otimes j})\subseteq\glob(X,\L^{\otimes(i+j)}(-D))\bigr\}
$$
using Lemma~\ref{division}.  We note that for this last step we can
use a small basepoint-free subspace of~$\glob(X,\L^{\otimes j})$
computed in advance.

\subsection{Decomposing divisors into prime divisors}

Let $X$ be a complete, smooth, geometrically connected curve of
genus~$g$ over a field~$k$, with a projective embedding via a line
bundle~$\L$ as in~\S\ref{representation-curve}.  The problem we are
now going to study is how to find the decomposition of a given divisor
on~$X$ as a linear combination of prime divisors.  We will see below
that this can be done if we are given the algebra~$S_X^{(h)}$ for
sufficiently large $h$ and if we are able to compute the primary
decomposition of a finite commutative $k$-algebra.  We have seen
in~\S\ref{primary-decomposition-radicals} that this is possible in the
case where $k$ is perfect and we have an algorithm for factoring
polynomials in one variable over~$k$.

Let $i$ be a positive integer, and let $D$ be an effective divisor
such that
$$
\deg D\le i\deg\L-2g+1.
$$
We view $D$ as a closed subscheme of~$X$ via the canonical closed
immersion
$$
j_D\colon D\to X.
$$
For every line bundle~${\cal M}$ on~$X$, the $k$-vector space
$\glob(D,j_D^*{\cal M})$ is in a natural way a free module of rank~one
over~$\glob(D,\O_D)$.  The multiplication map
$$
\mu_{i,i}\colon\glob(X,\L^{\otimes i})\times
\glob(X,\L^{\otimes i})\longrightarrow
\glob(X,\L^{\otimes 2i})
$$
descends to a bilinear map
$$
\mu_{i,i}^D\colon\glob(D,j_D^*\L^{\otimes i})\times
\glob(D,j_D^*\L^{\otimes i})\longrightarrow
\glob(D,j_D^*\L^{\otimes 2i})
$$
of free modules of rank~1 over~$\glob(D,\O_D)$.  This map is perfect
in the sense of~\S\ref{reconstructing-algebra-from-bilinear-map}.

We now assume that the graded $k$-algebra $S_X^{(h)}$ as
in~\S\ref{representation-curve} is given for some integer $h\ge2$.
From the subspace $\glob(X,\L^{\otimes i}(-D))$
of~$\glob(X,\L^{\otimes i})$ we can then determine
$\glob(D,j_D^*\L^{\otimes i})$ as a $k$-vector space by means of the
short exact sequence
$$
0\longrightarrow\glob(X,\L^{\otimes i}(-D))
\longrightarrow\glob(X,\L^{\otimes i})
\longrightarrow\glob(D,j_D^*\L^{\otimes i})
\longrightarrow0.
\eqnumber{exact-sequence-sections-on-divisor}
$$
(Note that exactness on the right follows from the assumption that
$\deg\L^{\otimes i}(-D)\ge 2g-1$.)  Similarly, we can compute
$\glob(D,j_D^*\L^{\otimes 2i})$ from $\glob(X,\L^{\otimes 2i}(-D))$
using the same sequence with $i$ replaced by $2i$.  We can then
determine the bilinear map~$\mu_{i,i}^D$ induced by~$\mu_{i,i}$ by
standard methods from linear algebra.

We then the method described
in~\S\ref{reconstructing-algebra-from-bilinear-map} to compute the
$k$-algebra~$\glob(D,\O_D)$ together with its action
on~$\glob(D,j_D^*\L^{\otimes i})$.  Next we determine the primary
decomposition of~$\glob(D,\O_D)$, say
$$
\glob(D,\O_D)\cong A_1\times A_2\times\cdots\times A_r,
$$
where each factor~$A_i$ is a finite local $k$-algebra with maximal
ideal~$P_i$; we assume the field~$k$ is such that we can do this (see
\S\ref{primary-decomposition-radicals}).  Such a prime ideal $P_i$
corresponds to a prime divisor in the support of~$D$, and the
corresponding multiplicity equals
$$
m_i={[A_i:k]\over[A_i/P_i:k]}.
$$

\algorithm (Decomposition of a divisor).  Let $X$ be a projective
curve over a field~$k$.  Let $i$ be a positive integer, and let $D$ be
an effective divisor such that
$$
\deg D\le i\deg\L_X-2\genus_X+1.
$$
Suppose that we have a (probabilistic) algorithm to compute the
primary decomposition of a finite commutative $k$-algebra~$A$ with
(expected) running time polynomial in~$[A:k]$, measured in operations
in~$k$.  Given the $k$-algebra~$S_X^{(2i)}$ and the
subspaces~$\glob(X,\L_X^{\otimes i}(-D))$ of $\glob(X,\L_X^{\otimes
i})$ and $\glob(X,\L_X^{\otimes 2i}(-D))$ of $\glob(X,\L_X^{\otimes
2i})$, this algorithm outputs the decomposition of~$D$ as a linear
combination of prime divisors as a list of pairs $(P,m_P)$, where $P$
is a prime divisor and $m_P$ is the multiplicity of~$P$ in~$D$.

\step Compute the spaces $\glob(D,j_D^*\L_X^{\otimes i})$ and
$\glob(D,j_D^*\L_X^{\otimes2i})$
using \eqref{exact-sequence-sections-on-divisor} and the analogous
short exact sequence with $2i$ in place of~$i$.

\step Compute the $k$-bilinear map $\mu_{i,i}^D$ from~$\mu_{i,i}$.

\step Compute a $k$-basis for $\glob(D,\O_D)$ as a linear subspace
of~$\End_k\glob(D,j_D^*\L_X^{\otimes i})$, where elements of the
latter $k$-algebra are expressed as matrices with respect to some
fixed basis of~$\glob(D,j_D^*\L_X^{\otimes i})$, as described
in~\S\ref{reconstructing-algebra-from-bilinear-map}.

\step Compute the multiplication table of~$\glob(D,\O_D)$ on the
$k$-basis of~$\glob(D,\O_D)$ found in the previous step.

\step Find the primary decomposition of~$\glob(D,\O_D)$.

\step For each local factor~$A$ computed in the previous step, let
$P_A$ denote the maximal ideal of~$A$, output the inverse image of
$P_A\cdot\glob(D,j_D^*\L_X^{\otimes i})$ in~$\glob(X,\L_X^{\otimes
i})$ and the integer $[A:k]\big/[A/P_A:k]$.

\endalgorithm

\label{decomposition-divisor}

\analysis It follows from the above discussion that the algorithm
returns the correct result.  It is straightforward to check that the
running time is polynomial in $i$ and~$\deg\L_X$, measured in
operations in~$k$.\endanalysis

A special case of this algorithm is when $D$ is the intersection
of~$X$ with a hypersurface of degree~$i-1$.  Let $s$ be a non-zero
section of~$\L_X^{\otimes(i-1)}$ defining this hypersurface.  The
subspaces that are used in this algorithm can then be computed as
$$
\eqalignno{
\glob(X,\L_X^{\otimes i}(-D))&=s\glob(X,\L_X)\cr
\noalign{\noindent and}
\glob(X,\L_X^{\otimes 2i}(-D))&=s\glob(X,\L_X^{\otimes(i+1)}).}
$$

\subsection{Finite morphisms between curves}

\label{finite-morphisms}

Let us now look at finite morphisms between curves.  A finite morphism
$$
f\colon X\to Y
$$
of complete, smooth, geometrically connected curves induces two
functors
$$
\eqalignno{
f^*\colon\{\hbox{line bundles on }Y\}&\to
\{\hbox{line bundles on }X\}\cr
\noalign{\noindent and}
\Norm_f\colon\{\hbox{line bundles on }X\}&\to
\{\hbox{line bundles on }Y\}.}
$$
Here $f^*{\cal N}$ denotes the usual inverse image of the line
bundle~${\cal N}$ on~$Y$, and $\Norm_f{\cal M}$ is the {\it norm\/} of
the line bundle~${\cal M}$ on~$X$ under the morphism~$f$.

Let us briefly explain the notion of the norm of a line bundle.  The
norm functor is a special case (that of ${\bf G}_{\rm m}$-torsors) of
the {\it trace of a torsor\/} under a finite locally free morphism;
see Deligne~\citex{SGA 4 III}{expos\'e~XVII, \nos 6.3.20--6.3.26}.  We
formulate the basic results for arbitrary finite locally free
morphisms of schemes
$$
f\colon X\to Y.
$$
In this situation there exists a functor
$$
\Norm_f\colon\{\hbox{line bundles on }X\}\to
\{\hbox{line bundles on }Y\}
$$
together with a collection of homomorphisms
$$
\Norm_f^\L\colon f_*\L\to\Norm_{E/k}\L
$$
of sheaves of sets, for all line bundles~$\L$ on~$X$, functorial under
isomorphisms of line bundles on~$X$, sending local generating sections
on~$X$ to local generating sections on~$Y$ and such that the equality
$$
\Norm_f^\L(xl)=\Norm_f(x)\cdot\Norm_f^\L(l)
$$
holds for all local sections $x$ of~$f_*\O_X$ and $l$ of $f_*\L$.
Here $\Norm_f\colon f_*\O_X\to\O_Y$ denotes the usual norm map for a
finite locally free morphism.  Moreover, the functor~$\Norm_f$
together with the collection of the $\Norm_f^\L$ is unique up to
unique isomorphism.  Instead of $\Norm_f$ we also write $\Norm_{X/Y}$
if the morphism~$f$ is clear from the context.

The basic properties of the norm functor are the following \citex{SGA
  4 III}{expos\'e~XVII, \no 6.3.26}:
\smallskip
\item{(1)} the functor $\Norm_f$ is compatible with any base
  change $Y'\to Y$;
\smallskip
\item{(2)} if $\L_1$ and~$\L_2$ are two line bundles on~$X$,
  there is a natural isomorphism
$$
\Norm_f(\L_1\otimes_{\O_X}\L_2)\cong
\Norm_f\L_1\otimes_{\O_Y}\Norm_f\L_2;
$$
\smallskip
\item{(3)} if $X\morphism f Y\morphism g Z$ are finite locally free
  morphisms, there is a natural isomorphism
$$
\Norm_{g\circ f}\isom\Norm_g\circ\Norm_f.
$$
\smallskip\noindent
Furthermore, there is a functorial isomorphism
$$
\Norm_f\L\isom\Hom_{\O_Y}(\det_{\O_Y}f_*\O_X,\det_{\O_Y}f_*\L);
\eqnumber{norm-det}
$$
see Deligne~\citex{SGA 4 III}{expos\'e~XVIII, \no 1.3.17}, and compare
Hartshorne~\citex{Hartshorne}{IV, Exercise~2.6}.

We now consider projective curves $X$ and~$Y$ as defined
in~\S\ref{representation-curve}.  Suppose we have a finite morphism
$$
f\colon X\to Y
$$
with the property that $f$ is induced by a graded homomorphism
$$
f^\#\colon S_Y\to S_X
$$
between the homogeneous coordinate rings of $Y$ and~$X$, or
equivalently by a morphism of the corresponding affine cones over $X$
and~$Y$.  Then $f^\#$ induces an isomorphism
$$
f^*\L_Y\isom\L_X
$$
of line bundles on~$X$; see Hartshorne~\citex{Hartshorne}{Chapter~II,
Proposition~5.12(c)}.  In particular, this implies
$$
\deg\L_X=\deg f\cdot\deg\L_Y.
$$

We represent a finite morphism~$f\colon X\to Y$ by the $k$-algebras
$S_X^{(h)}$ and~$S_Y^{(h)}$ for some $h\ge2$, together with the
$k$-algebra homomorphism
$$
f^\#\colon S_Y^{(h)}\to S_X^{(h)}
$$
induced by~$f^\#\colon S_Y\to S_X$, given as a collection of linear
maps $\glob(Y,\L_Y^{\otimes i})\to\glob(X,\L_X^{\otimes i})$
compatible with the multiplication maps on both sides.

In the following, when we mention a {\it finite morphism $f\colon X\to
Y$ between projective curves\/}, we assume that the $k$-algebras
$S_X^{(h)}$ and $S_Y^{(h)}$ and the homomorphism $f^\#\colon
S_Y^{(h)}\to S_X^{(h)}$ are given for some $h\ge2$.  A lower bound
for~$h$ will be specified in each of the algorithms that we describe.

\remark The homomorphism $f^\#$ gives rise to an injective $k$-linear
map
$$
\glob(Y,\L_Y)\to\glob(X,\L_X).
$$
Given this map, we can reconstruct $S(Y)$ as a subalgebra of~$S(X)$ by
noting that $S(Y)$ is generated as a $k$-algebra by~$\glob(Y,\L_Y)$.

\subsection{Images, pull-backs and push-forwards of divisors}

\label{functoriality-for-finite-morphisms}

Let us consider a finite morphism $f\colon X\to Y$ between complete,
smooth, geometrically connected curves over a field~$k$.  Such a
morphism~$f$ induces various maps between the groups of divisors
on~$X$ and on~$Y$.

First, for an {\it effective\/} divisor~$D$ on~$X$, we write $f(D)$
for the schematic image of~$D$ under~$f$.  The definition implies that
the ideal sheaf $\O_Y(-f(D))$ is the inverse image of~$f_*\O_X(-D)$
under the natural map $\O_Y\to f_*\O_X$.

Second, for any divisor~$D$ on~$X$, we have the ``push-forward''
$f_*D$ of~$D$ by~$f$; see Hartshorne~\citex{Hartshorne}{IV,
  Exercise~2.6}.  If $P$ is a prime divisor on~$X$, then its
image~$f(P)$ under~$f$ is a prime divisor on~$Y$, the residue
field~$k(P)$ is a finite extension of~$k(f(P))$, and $f_*P$ is given
by the formula
$$
f_*P=[k(P):k(f(P))]\cdot f(P).
\eqnumber{push-forward-prime-divisor}
$$
The residue field extension degree at~$P$ can simply be computed as
$$
\eqalign{
[k(P):k(f(P))]&={[k(P):k]\over[k(f(P)):k]}\cr
&={\deg P\over\deg f(P)}.}
$$

Third, for any divisor~$E$ on~$Y$, we have the ``pull-back'' $f^*E$
of~$E$ by~$f$; see for example
Hartshorne~\citex{Hartshorne}{page~137}.  If $Q$ is a prime divisor
on~$Y$, then $f^*Q$ is given by the formula
$$
f^*Q=\sum_{P\colon\,f(P)=Q}e(P)\cdot P
\eqnumber{pull-back-prime-divisor}
$$
where $P$ runs over the prime divisors of~$X$ mapping to~$Q$ and
$e(P)$ denotes the ramification index of~$f$ at~$P$.

Both $f_*$ and~$f^*$ are extended to arbitrary divisors on $X$ and~$Y$
by linearity.  Note that \eqref{push-forward-prime-divisor}
and~\eqref{pull-back-prime-divisor} imply the well-known formula
$$
f_*f^*E=(\deg f)E
$$
for any divisor~$E$ on~$Y$.  Furthermore, if $E$ is an {\it
effective\/} divisor on~$Y$, we have an equality
$$
f^*E=E\times_Y X
$$
of closed subschemes of~$X$, and if ${\cal I}_E$ denotes the ideal
sheaf of~$E$, then its inverse image $f^{-1}{\cal I}_E$ is the ideal
sheaf of~$f^*E$.

\remark The map $D\mapsto f(D)$ is not in general linear in~$D$.  We
do not extend it to the divisor {\it group\/} on~$X$, and in fact will
only need schematic images of {\it prime\/} divisors on~$X$ in what
follows.  In contrast, the maps $f_*$ and~$f^*$ are linear by
definition.

Now assume $f$ is a finite morphism between {\it projective\/} curves,
in the sense of~\S\ref{finite-morphisms}.  In particular, we have a
homomorphism $f^\#\colon S_Y\to S_X$ of graded~$k$-algebras.  We will
give algorithms to compute the image and the push-forward of a divisor
on~$X$ as well as the pull-back of a divisor on~$Y$.

The schematic image~$f(D)$ of an effective divisor~$D$ on~$X$ can be
computed using the following obvious algorithm.

\algorithm (Image of a divisor under a finite morphism).  Let $f\colon
X\to Y$ be a finite morphism between projective curves, let $i$ be a
positive integer, and let $D$ be an effective divisor on~$X$.  Given
the $k$-algebras $S_X^{(i)}$ and~$S_Y^{(i)}$, the homomorphism
$f^\#\colon S_Y^{(i)}\to S_X^{(i)}$ and the
subspace~$\glob(X,\L_X^{\otimes i}(-D))$ of~$\glob(X,\L_X^{\otimes
i})$, this algorithm outputs the subspace $\glob(Y,\L_Y^{\otimes
i}(-f(D)))$ of~$\glob(Y,\L_Y^{\otimes i})$.

\step Output the inverse image of the subspace~$\glob(X,\L_X^{\otimes
  i}(-D))$ of~$\glob(X,\L_X^{\otimes i})$ under the linear map
$\glob(Y,\L_Y^{\otimes i})\to\glob(X,\L_X^{\otimes i})$.

\endalgorithm

\label{image-divisor}

\analysis The definition of~$f(D)$ implies that $\L_Y^{\otimes
i}(-f(D))$ equals the inverse image of~$f_*\L_X^{\otimes i}(-D)$ under
the natural map $\L_Y^{\otimes i}\to f_*\L_X^{\otimes i}$.  Taking
global sections, we see that $\glob(Y,\L_Y^{\otimes i}(-f(D)))$ is the
inverse image of~$\glob(X,\L_X^{\otimes i}(-D))$ under the natural map
$\glob(Y,\L_Y^{\otimes i})\to\glob(X,\L_X^{\otimes i})$.  It is clear
that the algorithm needs a number of operations in~$k$ that is
polynomial in $\deg\L_X$ and~$i$.\endanalysis

\remark In the above algorithm, there are no restrictions on the
degrees of $D$ and~$f(D)$.  However, $f(D)$ is not uniquely determined
by $\glob(Y,\L_Y^{\otimes i}(-f(D)))$ if its degree is too large.

The algorithm to compute pull-backs that we will now give is based on
the fact that the pull-back of an effective divisor~$E$ is simply the
fibred product $E\times_Y X$, viewed as a closed subscheme of~$X$.  In
particular, the algorithm does not have to compute the ramification
indices, so instead we can {\it use\/} it to compute ramification
indices.  Namely, if $P$ is a prime divisor on~$X$, we see
from~\eqref{pull-back-prime-divisor} that the ramification index
at~$P$ equals the multiplicity with which $P$ occurs in the
divisor~$f^*(f(P))$.

\algorithm (Pull-back of a divisor under a finite morphism).  Let
$f\colon X\to Y$ be a finite morphism between projective curves.  Let
$i$ and~$j$ be positive integers, and let $E$ be an effective divisor
on~$Y$ such that
$$
\deg f\cdot\deg E\le i\deg\L_X-2\genus_X,\quad
\deg E\le i\deg\L_Y-2\genus_Y
$$
and
$$
(j-i)\deg\L_X+\deg f\cdot\deg E\ge 2\genus_X-1.
$$
(If we take $j\ge i+1$, the last equality does not pose an extra
restriction on~$E$.)  Given the $k$-algebras $S_X^{(i+j)}$
and~$S_Y^{(i+j)}$, the $k$-algebra homomorphism $f^\#\colon
S_Y^{(i+j)}\to S_X^{(i+j)}$ and the subspace $\glob(Y,\L_Y^{\otimes
i}(-E))$ of~$\glob(Y,\L_Y^{\otimes i})$, this algorithm outputs the
subspace $\glob(X,\L_X^{\otimes i}(-f^*E))$ of~$\glob(X,\L_X^{\otimes
i})$.

\step Compute the image $W$ of~$\glob(Y,\L_Y^{\otimes i}(-E))$ under
the linear map
$$
f^\#\colon\glob(Y,\L_Y^{\otimes i})\to\glob(X,\L_X^{\otimes i}).
$$

\step Compute the space $\glob(X,\L_X^{\otimes i+j}(-f^*E))$ as the
product of $W$ and $\glob(X,\L_X^{\otimes j})$ (see
Lemma \ref{lemma-mult-bpf}).

\step Compute $\glob(X,\L_X^{\otimes i}(-f^*E))$ using
Lemma~\ref{division}, and output the result.

\endalgorithm

\label{pull-back-divisor}

\analysis The ideal in $S_Y$ defining $E$ is generated by the linear
forms vanishing on~$E$, and the ideal of~$S_X$ defining $f^*E$ is
generated by the pull-backs of these forms.  This shows that $f^*E$ is
defined by the forms in~$W$.  In the second and third step, we compute
the space of all forms vanishing on~$f^*E$ is computed, i.e.\ the
inflation of~$W$.  That the method described is correct was proved
in~\S\ref{deflation-inflation}.  The running time is clearly
polynomial in $\deg\L_X$, $i$ and~$j$.\endanalysis

\algorithm (Push-forward of a divisor under a finite morphism).  Let
$f\colon X\to Y$ be a finite morphism between projective curves over a
field~$k$, let $i$ be a positive integer, and let $D$ be an effective
divisor on~$X$ such that
$$
\deg D\le i\deg\L_X-2\genus_X-1\quad\hbox{and}\quad
\deg D\le i\deg\L_Y-2\genus_Y.
$$
Suppose that we have a (probabilistic) algorithm to compute the
primary decomposition of a finite commutative $k$-algebra~$A$ with
(expected) running time polynomial in~$[A:k]$, measured in operations
in~$k$.  Given the $k$-algebras $S_X^{(2i)}$ and~$S_Y^{(2i)}$, the
homomorphism $f^\#\colon S_Y^{(2i)}\to S_X^{(2i)}$ and the
subspace~$\glob(X,\L_X^{\otimes i}(-D))$ of~$\glob(X,\L_X^{\otimes
i})$, this algorithm outputs the subspace~$\glob(Y,\L_Y^{\otimes
i}(-f_*D))$ of~$\glob(Y,\L_Y^{\otimes i})$.

\step Compute $\glob(X,\L_X^{\otimes 2i}(-D))$ as the product of
$\glob(X,\L_X^{\otimes i})$ and~$\glob(X,\L_X^{\otimes i}(-D))$ (see
Lemma~\ref{multiplication}).

\step Find the decomposition of~$D$ as a linear combination $\sum_P
n_PP$ of prime divisors using Algorithm~\ref{decomposition-divisor}.

\step For each prime divisor~$P$ in the support of~$D$,
compute the space $\glob(Y,\L^{\otimes i}(-f(P)))$ using
Algorithm~\ref{image-divisor}, and compute $[k(P):k(f(P))]$.

\step Compute the space $\glob(Y,\L_Y^{\otimes i}(-f_*D))$, where
$$
f_*D=\sum_P n_P[k(P):k(f(P))]f(P),
$$
and output the result.

\endalgorithm

\analysis The correctness of the algorithm follows from the definition
of~$f_*$.  It runs in (probabilistic) polynomial time in $\deg\L_X$
and~$i$, measured in field operations in~$k$.\endanalysis

We include here another algorithm that computes the push-forward of an
effective divisor under a non-constant rational function $X\to\P^1$ in
a slightly different setting than before.  We only assume $X$ to be
given as a projective curve, and we represent effective divisors
on~$\P^1$ as zero loci of homogeneous polynomials.  For simplicity, we
only consider divisors of degree at most $\deg\L_X$.

\algorithm (Push-forward of an effective divisor by a rational
function).  Let $X$ be a projective curve over a field~$k$, let $i$ be
a positive integer, let $\psi$ be a non-constant rational function
on~$X$ given as the quotient of two sections
$s,t\in\glob(X,\L_X^{\otimes i})$ without common zeroes, and let $D$
be an effective divisor on~$X$ of degree $d\le\deg\L_X$.  Suppose that
we have a (probabilistic) algorithm to compute the primary
decomposition of a finite commutative $k$-algebra~$A$ with (expected)
running time polynomial in~$[A:k]$, measured in operations in~$k$.
Given the $k$-algebra $S_X^{(\max\{4,i\})}$ and the subspace
$\glob(X,\L_X^{\otimes2}(-D))$, this algorithm outputs the homogeneous
polynomial of degree~$d$ defining the closed subscheme~$\psi_*D$
of~$\P^1_k$.  (This polynomial is unique up to multiplication by
elements of~$k^\times$)

\step Compute the space $\glob(X,\L_X^{\otimes4}(-D))$, and use
Algorithm~\ref{decomposition-divisor} to compute the decomposition
of~$D$ as a linear combination $D=\sum_Q n_Q Q$ of prime divisors.

\step For each prime divisor~$Q$ occurring in the decomposition
of~$D$:

\plusindent

\step Compute the base change $X_{k(Q)}$, where $k(Q)$ is the residue
field of~$Q$.  Compute the primary decomposition of~$Q_{k(Q)}$ and
pick a rational point $Q'$ in it.

\step Compute $\glob(X_{k(Q)},\L_X^{\otimes 2}(-Q'))$, then compute
the (one-dimensional) intersection of this space with $k\cdot s+k\cdot
t$, and express some generator of this intersection as $b_Qs-a_Qt$
with $a_Q,b_Q\in k(Q)$.  The element $\psi(Q')\in\P^1(k(Q))$ now has
homogeneous coordinates $(a_Q:b_Q)$.

\step Compute the homogeneous polynomial
$$
f_{\psi_*Q}=\Norm_{k(Q)/k}(b_Qu-a_Qv)\in k[u,v]
$$
defining $\psi_*Q$.

\minusindent

\step Output the homogeneous polynomial
$$
f_{\psi_*D}=\prod_Q f_{\psi_*Q}^{n_Q}\in k[u,v]
$$
of degree~$d$ defining $\psi_*D$.

\endalgorithm

\label{push-forward-rational-function}

\analysis It is straightforward to check that the algorithm is correct
and has expected running time polynomial in $i$ and $\deg\L_X$,
counted in operations in~$k$.\endanalysis

\subsection{The norm functor for effective divisors}

\label{norm-functor-effective-divisors}

Let $X$ be a proper, smooth, geometrically connected curve over a
field~$k$, and let $E$ be an effective divisor on~$X$.  We view $E$ as
a closed subscheme of~$X$, finite over~$k$, and we write
$$
j_E\colon E\to X
$$
for the closed immersion of~$E$ into~$X$.  For the purposes
of \S\ref{computing-frey-rueck-pairings} below, we will need an
explicit description of the norm functor $\Norm_{E/k}$ (for the
canonical morphism $E\to\Spec k$) that we saw
in~\S\ref{finite-morphisms}.  We view $\Norm_{E/k}$ as a functor from
free $\O_E$-modules of rank~1 to $k$-vector spaces of dimension~1.

Let ${\cal M}$ be a line bundle on~$X$.  We abbreviate
$$
\eqalignno{
\glob(E,{\cal M})&=\glob(E,j_E^*{\cal M})\cr
\noalign{\noindent and}
\Norm_{E/k}{\cal M}&=\Norm_{E/k}(j_E^*{\cal M}).}
$$
Suppose we have two line bundles ${\cal M}^+$ and~${\cal M}^-$, both
of degree at least $\deg E+2g-1$, together with an isomorphism
$$
{\cal M}\cong\shHom_{\O_X}({\cal M}^-,{\cal M}^+).
$$
Then we can compute $\glob(E,{\cal M}^-)$ and~$\glob(E,{\cal M}^+)$
using the short exact sequences
$$
0\longrightarrow\glob(X,{\cal M}^\pm(-E))
\longrightarrow\glob(X,{\cal M}^\pm)\longrightarrow
\glob(E,{\cal M}^\pm)\longrightarrow0,
$$
and we can express $\Norm_{E/k}$ via the isomorphism
$$
\Norm_{E/k}{\cal M}\cong\Hom_k\bigl(\det_k\glob(E,{\cal M}^-),
\det_k\glob(E,{\cal M}^+)\bigr)
\eqnumber{norm-det-pm}
$$
deduced from~\eqref{norm-det}.  We fix $k$-bases of $\glob(E,{\cal
M}^-)$ and $\glob(E,{\cal M}^+)$.  From the induced trivialisations of
$\det_k\glob(E,{\cal M}^\pm)$ we then obtain a trivialisation
of~$\Norm_{E/k}{\cal M}$.

Now consider three line bundles ${\cal M}$, ${\cal N}$ and~${\cal P}$,
together with an isomorphism
$$
\mu\colon{\cal M}\otimes_{\O_X}{\cal N}\isom{\cal P}.
$$
By the linearity of the norm functor, $\mu$ induces an isomorphism
$$
\Norm_{E/k}{\cal M}\otimes_k\Norm_{E/k}{\cal N}
\isom\Norm_{E/k}{\cal P}.
\eqnumber{lin-MNP}
$$
As above, we choose isomorphisms
$$
{\cal M}\cong\shHom_{\O_X}({\cal M}^-,{\cal M}^+),\quad
{\cal N}\cong\shHom_{\O_X}({\cal N}^-,{\cal N}^+),\quad
{\cal P}\cong\shHom_{\O_X}({\cal P}^-,{\cal P}^+)
$$
on $X$, where ${\cal M}^\pm$, ${\cal N}^\pm$ and ${\cal P}^\pm$ are
line bundles of degree at least $\deg E+2g+1$.  We fix bases of the
six $k$-vector spaces
$$
\glob(E,{\cal M}^\pm),\quad
\glob(E,{\cal N}^\pm),\quad
\glob(E,{\cal P}^\pm).
$$
Then \eqref{norm-det-pm} gives trivialisations of $\Norm_{E/k}{\cal
M}$, $\Norm_{E/k}{\cal N}$ and $\Norm_{E/k}{\cal P}$.  Under these
trivialisations, the isomorphism~\eqref{lin-MNP} equals multiplication
by some element $\lambda\in k^\times$.

To find an expression for~$\lambda$, we choose generators
$\alpha_{\cal M}^\pm$ and $\alpha_{\cal N}^\pm$ of the $\O_E$-modules
$\glob(E,{\cal M}^\pm)$ and~$\glob(E,{\cal N}^\pm)$.  To these we
associate the isomorphisms
$$
\eqalignno{
\alpha_{\cal M}\colon\glob(E,{\cal M}^-)&\isom\glob(E,{\cal M}^+)\cr
\noalign{\noindent and}
\alpha_{\cal N}\colon\glob(E,{\cal N}^-)&\isom\glob(E,{\cal N}^-)}
$$
sending $\alpha_{\cal M}^-$ to~$\alpha_{\cal M}^+$ and $\alpha_{\cal
N}^-$ to $\alpha_{\cal N}^+$, respectively.  Viewing $\alpha_{\cal M}$
and~$\alpha_{\cal N}$ as generators of~$\glob(E,{\cal M})$
and~$\glob(E,{\cal N})$ and applying the isomorphism
$$
\mu\colon\glob(E,{\cal M})\otimes_{\glob(E,{\O_E})}
\glob(E,{\cal N})\isom\glob(E,{\cal P})
$$
to~$\alpha_{\cal M}\otimes\alpha_{\cal N}$ we obtain a generator
of~$\glob(E,{\cal P})$, which we can identify with an isomorphism
$$
\alpha_{\cal P}\colon\glob(E,{\cal P}^-)\isom\glob(E,{\cal P}^+).
$$
We define $\delta_{\cal M}$ as the determinant of the matrix
of~$\alpha_{\cal M}$ with respect to the chosen bases.  Under the
given trivialisations of $\Norm_{E/k}{\cal M}$, the element
$\Norm_{E/k}^{\cal M}\alpha_{\cal M}$ corresponds to~$\delta_{\cal
M}$.  The same goes for ${\cal N}$ and~${\cal P}$.  On the other hand,
the isomorphism~\eqref{lin-MNP} maps $\Norm_{E/k}^{\cal M}\alpha_{\cal
M}\otimes\Norm_{E/k}^{\cal N}\alpha_{\cal N}$ to~$\Norm_{E/k}^{\cal
P}\alpha_{\cal P}$.  We conclude that we can express $\lambda$ as
$$
\lambda={\delta_{\cal P}\over\delta_{\cal M}\delta_{\cal N}}.
\eqnumber{lambda-deltadeltadelta}
$$

Let us turn the above discussion into an algorithm.  Let $X$ be a
projective curve over~$k$, embedded via a line bundle~$\L$, and let
$E$ be an effective divisor on~$X$.  For simplicity, we restrict to
the case where
$$
\deg E\le\deg\L.
$$
We consider line bundles
$$
{\cal M}=\L^{\otimes i}(-D_1)\quad\hbox{and}\quad
{\cal N}=\L^{\otimes j}(-D_2),
$$
where $i$ and $j$ are non-negative integers and $D_1$ and~$D_2$ are
effective divisors such that
$$
\deg D_1=i\deg\L\quad\hbox{and}\quad\deg D_2=j\deg\L.
$$
We take
$$
{\cal M}^-={\cal N}^-={\cal P}^-=\L^{\otimes2}
$$
and
$$
\displaylines{
{\cal M}^+=\L^{\otimes(i+2)}(-D_1),\quad
{\cal N}^+=\L^{\otimes(j+2)}(-D_2),\cr
{\cal P}^+=\L^{\otimes(i+j+2)}(-D_1-D_2).}
$$

\algorithm (Linearity of the norm functor). Let $X$ be a projective
curve over a field~$k$, and let $E$, $D_1$ and~$D_2$ be effective
divisors on~$X$ such that
$$
\deg E=\deg\L,\quad
\deg D_1\le i\deg\L,\quad
\deg D_2\le j\deg\L.
$$
Fix bases of the four $k$-vector spaces
$$
\displaylines{
\glob(E,\L^{\otimes2}),\quad
\glob(E,\L^{\otimes(i+2)}(-D_1)),\cr
\glob(E,\L^{\otimes(j+2)}(-D_2)),\quad
\glob(E,\L^{\otimes(i+j+2)}(-D_1-D_2)).}
$$
and consider the corresponding trivialisations
$$
\displaylines{
t_1\colon k\isom \Norm_{E/k}\L^{\otimes i}(-D_1),\quad
t_2\colon k\isom \Norm_{E/k}\L^{\otimes j}(-D_2),\cr
t_3\colon k\isom \Norm_{E/k}\L^{\otimes i+j}(-D_1-D_2)}
$$
defined by~\eqref{norm-det-pm}.  Given the $k$-algebra
$S_X^{(i+j+4)}$, bases for the $k$-vector spaces
$$
\displaylines{
\glob(X,\L^{\otimes2}),\quad
\glob(X,\L^{\otimes(i+2)}),\cr
\glob(X,\L^{\otimes(j+2)}(-D_2)),\quad
\glob(X,\L^{\otimes(i+j+2)}(-D_1-D_2))}
$$
and the quotient maps
$$
\eqalign{
\glob(X,\L^{\otimes2})&\longrightarrow\glob(E,\L^{\otimes2}),\cr
\glob(X,\L^{\otimes(i+2)}(-D_1))&\longrightarrow
  \glob(E,\L^{\otimes i+2}(-D_1)),\cr
\glob(X,\L^{\otimes(j+2)}(-D_2))&\longrightarrow
  \glob(E,\L^{\otimes j+2}(-D_2)),\cr
\glob(X,\L^{\otimes(i+j+2)}(-D_1-D_2))&\longrightarrow
  \glob(E,\L^{\otimes i+2}(-D_1))\cr}
$$
as matrices with respect to the given bases, this algorithm outputs
the element $\lambda\in k^\times$ such that the diagram
$$
\commdiag{
k& \isomorphism{t_1\otimes t_2}& \Norm_{E/k}\L^{\otimes i}(-D_1)
\otimes_k\Norm_{E/k}\L^{\otimes j}(-D_2)\cr
\leftlabel\lambda\downar\rightlabel\sim& & \downar\rightlabel\sim\cr
k& \isomorphism{t_3}& \Norm_{E/k}\L^{\otimes(i+j)}(-D_1-D_2)\cr}
$$
is commutative.

\step Compute the spaces
$$
\glob(E,\L^{\otimes(i+4)}(-D_1))\quad\hbox{and}\quad
\glob(E,\L^{\otimes(i+j+4)}(-D_1-D_2))
$$
and the multiplication maps
$$
\displaylines{
\glob(E,\L^{\otimes2})\times\glob(E,\L^{\otimes(i+2)}(-D_1))\to
\glob(E,\L^{\otimes(i+4)}(-D_1)),\cr
\glob(E,\L^{\otimes(i+2)}(-D_1))\times\glob(E,\L^{\otimes(j+2)}(-D_2))\to
\glob(E,\L^{\otimes(i+j+4)}(-D_1-D_2)),\cr
\glob(E,\L^{\otimes2})\times\glob(E,\L^{\otimes(i+j+2)}(-D_1-D_2))\to
\glob(E,\L^{\otimes(i+j+4)}(-D_1-D_2)).}
$$

\step Apply the probabilistic method described
in~\S\ref{reconstructing-algebra-from-bilinear-map} to the bilinear
maps just computed to find generators $\beta_0$, $\beta_1$
and~$\beta_2$ of the free $\glob(E,\O_E)$-modules
$\glob(E,\L^{\otimes2})$, $\glob(E,\L^{\otimes(i+2)}(-D_1))$ and
$\glob(E,\L^{\otimes(j+2)}(-D_2))$ of rank~1.
\hfill\break
(Note that we do not need the $k$-algebra structure
on~$\glob(E,\L^{\otimes2})$.  If $k$ is small, we may have to extend
the base field, but it is easy to see that this is not a problem.)

\step Compute the matrix (with respect to the given bases) of the
isomorphism~$\alpha_1$ defined by the commutative diagram
$$
\commdiag{
\glob(E,\L^{\otimes2})& \isomorphism{\alpha_1}&
  \glob(E,\L^{\otimes(i+2)}(-D_1))\cr
\bigm\Vert& & \leftlabel{\sim}\downar\rightlabel{\cdot\beta_0}\cr
\glob(E,\L^{\otimes2})& \isomorphism{\cdot\beta_1}&
  \glob(E,\L^{\otimes(i+4)}(-D_1))\rlap,}
$$
of the isomorphism~$\alpha_2$ defined by the similar diagram for
$\L^{\otimes j}(-D_2)$ instead of $\L^{\otimes i}(-D_1)$ and of the
isomorphism~$\alpha_3$ defined by the commutative diagram
$$
\commdiag{
\glob(E,\L^{\otimes2})& \isomorphism{\alpha_3}&
  \glob(E,\L^{\otimes(i+j+2)}(-D_1-D_2))\cr
\leftlabel{\alpha_1}\downar\rightlabel{\sim}& &
\leftlabel{\sim}\downar\rightlabel{\cdot\beta_0}\cr
  \glob(E,\L^{\otimes(i+2)}(-D_1))& \isomorphism{\cdot\beta_2}&
  \glob(E,\L^{\otimes(i+j+4)}(-D_1-D_2))\rlap.}
$$

\step Compute the elements $\delta_1$, $\delta_2$ and~$\delta_3$
of~$k^\times$ as the determinants of the matrices of $\alpha_1$,
$\alpha_2$ and~$\alpha_3$ computed in the previous step.

\step Output the element
$\displaystyle{\delta_3\over\delta_1\delta_2}\in k^\times$.

\endalgorithm

\label{algorithm-linearity-norm}

\analysis We note that $\beta_0$ plays the role of $\alpha_{\cal
M}^-$, $\alpha_{\cal N}^-$ and $\alpha_{\cal P}^-$ in the notation of
the discussion preceding the algorithm, and that $\beta_1$, $\beta_2$
and $\beta_1\beta_2/\beta_0$ play the roles of $\alpha_{\cal M}^+$,
$\alpha_{\cal N}^+$ and $\alpha_{\cal P}^+$.  This means that
$\alpha_1$, $\alpha_2$ and~$\alpha_3$ are equal to $\alpha_{\cal M}$,
$\alpha_{\cal N}$ and~$\alpha_{\cal P}$.  It now follows
from~\eqref{lambda-deltadeltadelta} that the output of the algorithm
is indeed equal to~$\lambda$.  It is clear that the algorithm runs in
(probabilistic) polynomial time in $\deg\L$, $i$ and~$j$, measured in
field operations in~$k$.\endanalysis

\subsection{Computing in the Picard group of a curve}

\label{computing-in-picard-group}

We now explain how to compute with elements in the Picard group of a
curve~$X$, using the operations on divisors described in the first
part of this section.  We only consider the group~$\Pic^0 X$ of
isomorphism classes of line bundles of degree~0.  This group can be
identified in a canonical way with a subgroup of rational points of
the Jacobian variety of~$X$.  If $X$ has a rational point, then this
subgroup consists of all the rational points of the Jacobian.

We will only describe Khuri-Makdisi's {\it medium model\/} of~$\Pic^0
X$ relative to a fixed line bundle~$\L$ of degree
$$
\deg\L\ge 2g+1,
$$
but at the same time
$$
\deg\L\le c(g+1)
$$
for some constant $c\ge 1$, as described in
Khuri-Makdisi~\citex{Khuri-Makdisi: Linear algebra algorithms}{\S5}.

\remark Khuri-Makdisi starts with a divisor $D_0$ whose degree
satisfies the above inequalities and takes $\L=\O_X(D_0)$.  This is of
course only a matter of language.  Another difference in notation is
that Khuri-Makdisi writes $\L_0$ for~$\L$ and uses the notation~$\L$
for~$\L_0^{\otimes2}$ (in the medium model) or~$\L_0^{\otimes3}$ (in
the {\it large\/} and {\it small\/} models, which we do not describe
here).

We represent elements of~$\Pic^0 X$ by effective divisors of
degree~$\deg\L$ as follows: the isomorphism class of a line
bundle~${\cal M}$ of degree~0 is represented by the divisor of some
global section of the line bundle ${\cal H\mit om}({\cal M},\L)$ of
degree~$\deg\L$, i.e.\ by any effective divisor~$D$ such that
$$
{\cal M}\cong\L(-D).
$$
It follows from the inequality $\deg\L\ge 2g$ that we can represent
any effective divisor~$D$ of degree~$\deg\L$ by the
subspace~$\glob(X,\L^{\otimes2}(-D))$ of codimension~$\deg\L$
in~$\glob(X,\L^{\otimes2})$.

There are a few basic operations:
\smallskip
\item{$\bullet$} {\it membership test\/}: given a subspace of
codimension~$\deg\L$ in~$\glob(X,\L^{\otimes2})$, decide whether it
represents an element of~$\Pic^0 X$, i.e.\ whether it is of the form
$\glob(X,\L^{\otimes2}(-D))$ for an effective divisor~$D$ of
degree~$\deg\L$.
\smallskip
\item{$\bullet$} {\it zero test\/}: given a subspace of
codimension~$\deg\L$ in~$\glob(X,\L^{\otimes2})$, decide whether it
represents the zero element of~$\Pic^0 X$.
\smallskip
\item{$\bullet$} {\it zero element\/}: output a subspace of
codimension~$\deg\L$ in~$\glob(X,\L^{\otimes2})$ representing the
element $0\in\Pic^0 X$.
\smallskip
\item{$\bullet$} {\it addflip\/}: given two subspaces
of~$\glob(X,\L^{\otimes2})$ representing elements $x,y\in\Pic^0 X$,
compute a subspace of~$\glob(X,\L^{\otimes2})$ representing the
element~$-x-y$.
\smallskip
\noindent From the ``addflip'' operation, one immediately gets
negation ($-x=-x-0$), addition ($x+y=-(-x-y)$) and subtraction
($x-y=-(-x)-y$).  Clearly, one can test whether two elements $x$
and~$y$ are equal by computing $x-y$ and testing whether the result
equals zero.

\remark With regard to actual implementations of the above algorithms,
we note that some of the operations can be implemented in a more
efficient way than by composing the basic operations just described.
We refer to~\cite{Khuri-Makdisi} for details.

By Khuri-Makdisi's results in~\cite{Khuri-Makdisi}, the above
operations can be implemented using randomised algorithms with
expected running time of~$O(g^{3+\epsilon})$ for any $\epsilon>0$,
measured in operations in the field~$k$.  This can be improved to
$O(g^{2.376})$ by means of fast linear algebra algorithms.  (The
exponent~2.376 is an upper bound for the complexity of matrix
multiplication.)

Multiplication by an integer~$n$ can be done efficiently by means of
an {\it addition chain\/} for $n$.  This is a sequence of positive
integers $(a_1,a_2,\ldots,a_m)$ with $a_1=1$ and $a_m=n$ such that for
each $l>1$ there exist $i(l)$ and~$j(l)$ in $\{1,2,\ldots,l-1\}$ such
that $a_l=a_{i(l)}+a_{j(l)}$.  We consider the indices $i(l)$
and~$j(l)$ as given together with the addition chain.  The
integer~$m$ is called the {\it length\/} of the addition chain.  A
more general and often slightly more efficient method of multiplying
by~$n$ is to use an {\it addition-subtraction chain\/}, where $a_l$ is
allowed to be either $a_{i(l)}+a_{j(l)}$ or $a_{i(l)}-a_{j(l)}$.
However, since the ``addflip'' operation in our set-up takes less time
than the addition or subtraction algorithms, the most worthwhile
option is to use an {\it anti-addition chain\/}, which is a sequence
of (not necessarily positive) integers $(a_0,a_1,\ldots,a_m)$ such
that
$$
a_l=\cases{0& if $l=0$;\cr
1& if $l=1$;\cr
-a_{i(l)}-a_{j(l)}& if $2\le l\le m$}
$$
and $a_m=n$; the $i(l)$ and~$j(l)$ are given elements
of~$\{0,1,\ldots,l-1\}$ for $2\le l\le m$.

It is well known that for every positive integer~$n$ there exists an
addition chain whose length is bounded by a constant times~$\log n$.
Moreover, there are algorithms (such as the {\it binary method\/} used
for repeated squaring) to find such an addition chain in time $O((\log
n)^2)$.  We leave it to the reader to write down a similar algorithm
for finding an anti-addition chain.

For later use, we give versions of the ``zero test'' and ``addflip''
algorithms that are identical to those given by Khuri-Makdisi, except
that some extra information computed in the course of the algorithm is
part of the output.

\algorithm (Zero test).  Let $X$ be a projective curve over a
field~$k$, and let $x$ be an element of~$\Pic^0 X$.  Given the
$k$-algebra $S_X^{(2)}$ and a subspace $\gl(X,\L_X^{\otimes2}(-D))$
of~$\gl(X,\L_X^{\otimes2})$ representing $x$, this algorithm outputs
{\tt false} if $x\ne0$ (i.e.\ if the line bundle~$\L_X(-D)$ is
non-trivial).  If $\L_X(-D)$ is trivial, the algorithm outputs a pair
$({\tt true},s)$, where $s$ is a global section of~$\L_X$ with
divisor~$D$.

\label{zero-test}

\step Compute the space
$$
\gl(X,\L_X(-D))=\bigl\{s\in\gl(X,\L_X)\bigm|
s\gl(X,\L_X)\subseteq\glob(X,\L_X^{\otimes2}(-D))\bigr\}.
$$
(The truth of this equality follows from Lemma~\ref{division}.)

\step If $\gl(X,\L_X(-D))=0$, output {\tt false}.  Otherwise, output
$({\tt true},s)$, where $s$ is any non-zero element of the
one-dimensional $k$-vector space $\gl(X,\L_X(-D))$.

\endalgorithm

\algorithm (Addflip).  Let $X$ be a projective curve over a field~$k$,
and let $x$ and~$y$ be elements of~$\Pic^0 X$.  Given the $k$-algebra
$S_X^{(5)}$ and subspaces $\gl(X,\L_X^{\otimes2}(-D))$ and
$\gl(X,\L_X^{\otimes2}(-E))$ of~$\gl(X,\L_X^{\otimes2})$ representing
$x$ and~$y$, this algorithm outputs a subspace
$\gl(X,\L_X^{\otimes2}(-F))$ representing $-x-y$, as well as a global
section~$s$ of~$\L_X^{\otimes3}$ such that
$$
\div s=D+E+F.
$$

\label{addflip}

\step Compute $\gl(X,\L_X^{\otimes4}(-D-E))$ as the product of
$\gl(X,\L_X^{\otimes2}(-D))$ and~$\gl(X,\L_X^{\otimes2}(-E))$ (see
Lemma~\ref{multiplication}).

\step Compute the space
$$
\gl(X,\L_X^{\otimes3}(-D-E))=\bigl\{
s\in\gl(X,\L_X^{\otimes3})\bigm|
s\gl(X,\L_X)\subseteq\gl(X,\L_X^{\otimes4}(-D-E))\bigr\}
$$
(see Lemma~\ref{division}).

\step Choose any non-zero $s\in\gl(X,\L_X^{\otimes3}(-D-E))$.  Let $F$
denote the divisor of~$s$ as a global section
of~$\L_X^{\otimes3}(-D-E)$.

\step Compute the space
$$
\gl(X,\L_X^{\otimes5}(-D-E-F))=s\gl(X,\L_X^{\otimes2}).
$$

\step Compute the space
$$
\eqalign{
\gl(X,\L_X^{\otimes2}(-F))&=\bigl\{
t\in\gl(X,\L_X^{\otimes2})\bigm|\cr
&\qquad t\gl(X,\L_X^{\otimes3}(-D-E))\subseteq
\gl(X,\L_X^{\otimes5}(-D-E-F))\bigr\}}
$$
(see again Lemma~\ref{division}).

\step Output the space $\gl(X,\L_X^{\otimes2}(-F))$ and the
section~$s\in\gl(X,\L_X^{\otimes3})$.

\endalgorithm

\subsection{Normalised representatives of elements of the Picard group}

\label{normalised-representatives}

Let $X$ be a projective curve over a field~$k$, and let $O$ be a
$k$-rational point of~$X$.  Let $x$ be an element of~$\Pic^0 X$, and
let ${\cal M}$ be a line bundle representing $x$.  Let $r_x^{\L_X,O}$
be the greatest integer~$r$ such that
$$
\gl(X,\shHom({\cal M},\L_X(-rO)))\ne0.
$$
Then $\gl(X,\shHom_{\O_X}({\cal M},\L_X(-r_x^{\L_X,O}O)))$ is
one-dimensional, so there exists a unique effective divisor~$R$ such
that
$$
{\cal M}\cong\L_X(-R-r_x^{\L_X,O}O).
$$
We define the {\it $(\L_X,O)$-normalised representative\/} of~$x$ as
the effective divisor
$$
R^{\L_X,\O}_x=R+r_x^{\L_X,O}O
$$
of degree~$\deg\L_X$; it is a canonically defined divisor (depending
on~$O$) with the property that $x$ is represented
by~$\L_X(-R^{\L_X,O}_x)$.

\remark Since for any line bundle~${\cal N}$ we have
$$
\eqalignno{
\deg{\cal N}\ge g\;&\Longrightarrow\;\gl(X,{\cal N})\ne0\cr
\noalign{\noindent and}
\deg{\cal N}<0\;&\Longrightarrow\;\gl(X,{\cal N})=0,}
$$
the integer~$r_x^{\L_X,O}$ satisfies
$$
\deg\L_X-\genus_X\le r_x^{\L_X,O}\le\deg\L_X.
$$

\algorithm (Normalised representative). Let $X$ be a projective curve
over a field~$k$, and let $O$ be a $k$-rational point of~$X$.  Let $x$
be an element of~$\Pic^0 X$, and let $R^{\L_X,O}_x$ be the
$(\L_X,O)$-normalised representative of~$x$.  Given the $k$-algebra
$S_X^{(4)}$, the space $\gl(X,\L_X^{\otimes2}(-O))$ and a subspace
of~$\gl(X,\L_X^{\otimes2})$ representing $x$, this algorithm outputs
the integer $r^{\L_X,O}_x$ and the subspace
$\gl(X,\L_X^{\otimes2}(-R^{\L_X,O}_x))$ of~$\gl(X,\L_X^{\otimes2})$.

\step Using the negation algorithm, find a
subspace~$\gl(X,\L_X^{\otimes2}(-D))$ of~$\gl(X,\L_X^{\otimes2})$
representing $-x$.  Put $r=\deg\L_X$.

\step Compute $\gl(X,\L_X^{\otimes2}(-rO))$, then compute
$\gl(X,\L_X^{\otimes4}(-D-rO))$ as the product of
$\gl(X,\L_X^{\otimes2}(-D))$ and~$\gl(X,\L_X^{\otimes2}(-rO))$, and
then compute
$$
\gl(X,\L_X^{\otimes2}(-D-rO))=\bigl\{
t\in\gl(X,\L_X^{\otimes2})\bigm|
t\gl(X,\L_X^{\otimes2})\subseteq
\gl(X,\L_X^{\otimes4}(-D-rO))\bigr\}.
$$

\label{normalised-representative-iteration}

\step If $\gl(X,\L_X^{\otimes2}(-D-rO))=0$, decrease $r$ by~1 and go
to step~\ref{normalised-representative-iteration}.

\step Let $s$ be a non-zero element
of~$\gl(X,\L_X^{\otimes2}(-D-rO))$.  Compute
$$
\gl(X,\L_X^{\otimes4}(-D-R^{\L_X,O}_x))=s\gl(X,\L_X^{\otimes2}),
$$
and then compute
$$
\gl(X,\L_X^{\otimes2}(-R^{\L_X,O}_x))=\bigl\{
t\in\gl(X,\L_X^{\otimes2})\bigm|t\gl(X,\L_X^{\otimes2}(-D))\subseteq
\gl(X,\L_X^{\otimes4}(-D-R^{\L_X,O}_x))\bigr\},
$$

\step Output $r_x^{\L_X,O}=r$ and $\gl(X,\L_X^{\otimes2}(-R^{\L_X,O}_x))$.

\endalgorithm

\analysis It follows from the definition of~$R^{\L_X,O}_x$ that this
algorithm is correct.  It is straightforward to check that its running
time, measured in operations in~$k$, is polynomial
in~$\deg\L_X$.\endanalysis


\subsection{Descent of elements of the Picard group}

Now let $k'$ be a finite extension of~$k$, and write
$$
X'=X\times_{\Spec k}\Spec k'.
$$
Consider the natural inclusion map
$$
i\colon\Pic^0 X\to\Pic^0 X'.
$$
Let $x'$ be an element of~$\Pic^0 X'$.  We can use normalised
representatives to decide whether $x'$ lies in the image of $i$, and
if so, to find the unique element $x\in\Pic^0 X$ such that $x'=i(x)$.

\algorithm (Descent). Let $X$ be a projective curve over a field~$k$,
and let $O$ be a $k$-rational point of~$X$.  Let $k'$ be a finite
extension of~$k$, write
$$
X'=X\times_{\Spec k}\Spec k',
$$
and let $\L_{X'}$ denote the pull-back of the line bundle~$\L_X$
to~$X'$.  Let $x'$ be an element of~$\Pic^0 X'$.  Given the
$k$-algebra $S_X^{(4)}$, the spaces
$$
\glob(X,\L_X^{\otimes2}(-rO))\quad\hbox{for }
\deg\L_X-\genus_X\le d\le\deg\L_X
$$
and a subspace of~$\glob(X',\L_{X'}^{\otimes2})$ representing $x'$,
this algorithm outputs {\tt false} if $x'$ is in not the image of the
canonical map
$$
i\colon\Pic^0 X\to\Pic^0 X'.
$$
Otherwise, the algorithm outputs $({\tt
true},\glob(X,\L_X^{\otimes2}(-D)))$, where
$\glob(X,\L_X^{\otimes2}(-D))$ represents the unique element
$x\in\Pic^0 X$ such that $i(x)=x'$.

\step Compute the $(\L_X,O)$-normalised
representative~$R_{x'}^{\L_X,O}$ of~$x'$.

\step Compute the $k$-vector space
$$
V=\glob(X',\L_{X'}^{\otimes2}(-R_x))\cap\glob(X,\L_X^{\otimes2}).
$$

\step If the codimension of~$V$ in~$\glob(X,\L_X^{\otimes2})$ is
less that~$\deg\L_X$, output {\tt false}; otherwise, output $({\tt
true},V)$.

\endalgorithm

\label{descent}

\analysis In step~3, we check whether $R^{\L_X,O}_x$ is defined
over~$k$ or, equivalently, whether $x$ is defined over~$k$.  If this
is the case, the space~$V$ equals $\glob(X,\L_X^{\otimes2}(-R_x))$,
where $x$ is the unique element of~$\Pic^0 X$ such that $i(x)=x'$.
This shows that the algorithm is correct; its running time, measured
in operations in $k$ and~$k'$, is clearly polynomial in
$\deg\L_X$.\endanalysis

\subsection{Picard and Albanese maps}

\label{picard-albanese}

A finite morphism
$$
f\colon X\to Y
$$
between complete, smooth, geometrically connected curves over a
field~$k$ induces two group homomorphisms
$$
\eqalignno{
\Pic f\colon\Pic^0 Y&\to\Pic^0 X\cr
\noalign{\noindent and}
\Alb f\colon\Pic^0 X&\to\Pic^0 Y,}
$$
called the {\it Picard\/} and {\it Albanese\/} maps, respectively.  In
terms of line bundles, they can be described as follows.  The Picard
map sends the class of a line bundle~${\cal N}$ on~$Y$ to the class of
the line bundle~$f^*{\cal N}$ on~$X$, and the Albanese map sends the
class of a line bundle~${\cal M}$ on~$X$ to the class of the line
bundle~$\Norm_f{\cal M}$ on~$Y$.

Alternatively, these maps can be described in terms of divisor classes
as follows.  The group homomorphisms
$$
f_*\colon\Div^0 X\to\Div^0 Y\quad\hbox{and}\quad
f^*\colon\Div^0 Y\to\Div^0 X
$$
between the groups of divisors of degree~0 on $X$ and~$Y$ respect the
relation of linear equivalence on both sides.  The Picard map sends
the class of a divisor~$E$ on~$Y$ to the class of the divisor~$f^*E$
on~$X$, and the Albanese map sends the class of a divisor~$D$ on~$X$
to the class of the divisor~$f_*D$ on~$Y$.

Let us now assume that $f\colon X\to Y$ is a finite morphism of {\it
projective\/} curves in the sense of~\S\ref{finite-morphisms}.  The
following algorithms can be used to compute the maps $\Pic f$
and~$\Alb f$.  The algorithm for the Albanese map is mostly a wax
nose, since we only reduce the problem to a different one, namely that
of computing traces in Picard groups with respect to finite extensions
of the base field.  However, this is a problem that can be solved at
least for finite fields, as we will see
in~\S\ref{frobenius-endomorphism}.

\algorithm (Picard map). Let $f\colon X\to Y$ be a finite morphism of
projective curves, and let $y$ be an element of~$\Pic^0 Y$.  Given the
$k$-algebras $S_X^{(4)}$ and~$S_Y^{(4)}$, the homomorphism $f^\#\colon
S_Y^{(4)}\to S_X^{(4)}$ and a subspace $\glob(Y,\L_Y^{\otimes2}(-E))$
of~$\glob(Y,\L_Y^{\otimes2})$ representing $y$, this algorithm outputs
a subspace of~$\glob(X,\L_X^{\otimes2})$ representing $(\Pic
f)(y)\in\Pic^0 X$.

\step Compute the subspace~$\glob(X,\L_X^{\otimes2}(-D))$ for the
divisor $D=f^*E$ using Algorithm~\ref{pull-back-divisor} (taking
$i=j=2$ in the notation of that algorithm), and output the result.

\endalgorithm

\analysis Since $(\Pic f)(y)$ is represented by the line
bundle~$\L_X(-f^*D)$, the correctness of this algorithm follows from
that of Algorithm~\ref{pull-back-divisor}.  Furthermore, the running
time of Algorithm~\ref{pull-back-divisor}, measured in operations
in~$k$, is polynomial in $\deg\L_X$ for fixed $i$ and~$j$; therefore,
the running time of this algorithm is also polynomial
in~$\deg\L_X$.\endanalysis

\algorithm (Albanese map). Let $f\colon X\to Y$ be a finite morphism
of projective curves over a field~$k$.  Let $x$ be an element
of~$\Pic^0 X$, and let $O$ be a $k$-rational point of~$Y$.  Suppose
that we have a (probabilistic) algorithm to compute the primary
decomposition of a finite commutative $k$-algebra~$A$ with (expected)
running time polynomial in~$[A:k]$, measured in operations in~$k$.
Suppose furthermore that we can compute the trace of an element
$y\in\Pic^0(Y_{k'})$ over~$k$ for a finite extension $k'$ of~$k$ in
time polynomial in $\deg\L_Y$ and~$[k':k]$, measured in operations
in~$k$.  Given the $k$-algebras $S_X^{(6)}$ and~$S_Y^{(6)}$, the
homomorphism $f^\#\colon S_Y^{(6)}\to S_X^{(6)}$, the
space~$\glob(Y,\L_Y^{\otimes2}(-O))$ and a
subspace~$\glob(X,\L_X^{\otimes2}(-D))$ of~$\glob(Y,\L_Y^{\otimes2})$
representing $x$, this algorithm outputs a subspace
of~$\glob(Y,\L_Y^{\otimes2})$ representing $(\Alb f)(x)\in\Pic^0 Y$.

\step Compute $\glob(X,\L_X^{\otimes4}(-D))$ as the product
of $\glob(X,\L_X^{\otimes2})$ and $\glob(X,\L_X^{\otimes2}(-D))$.

\step Find the decomposition of~$D$ as a linear combination $\sum_P
n_PP$ of prime divisors using Algorithm~\ref{decomposition-divisor}.

\step For each $P$ occurring in the support of~$D$:

\plusindent

\step Compute the base changes $X_{k(P)}$ and~$Y_{k(P)}$.

\step Find the primary decomposition of the divisor $P_{k(P)}$
on~$X_{k(P)}$, and pick a rational point~$P'$ in it.

\step Compute the space
$\glob(Y_{k(P)},\L_Y^{\otimes2}(-f(P')-(\deg\L_Y-1)O))$; this
represents an element $y_{P'}\in\Pic^0(Y_{k(P)})$.

\step Compute the element $y_P=\tr_{k(P)/k}y_{P'}$ of~$\Pic^0
Y_{k(P)}$.  Apply Algorithm~\ref{descent} to get a representation
for~$y_P$ as an element of~$\Pic^0 Y$.

\minusindent

\step Compute the element $y=\sum_P n_Py_P$ of~$\Pic^0(Y)$.

\step Output the element $y-(\deg f)(\deg\L_Y-1)y_0$ of~$\Pic^0 Y$,
where $y_0$ is the element of~$\Pic^0 Y$ represented by
$\glob(Y,\L_Y^{\otimes2}(-(\deg\L_Y)O))$.

\endalgorithm

\analysis The definition of~$y_{P,i}$ implies that
$$
y_{P'}=[\L_Y(-f(P')-(\deg\L_Y-1)O)],
$$
the definition of~$y_P$ implies that
$$
y_P=[\L_Y^{\otimes[k(P):k]}(-f_*P-[k(P):k](\deg\L_Y-1)O)]
$$
and the definition of~$y$ implies that
$$
\eqalign{
y&=[\L_Y^{\otimes\deg\L_X}(-f_*D-(\deg\L_X)(\deg\L_Y-1)O)]\cr
&=[\L_Y^{\deg f}(-f_*D)]+(\deg f)(\deg\L_Y-1)[\L_Y(-(\deg\L_Y)O)].}
$$
Together with the definition of~$y_0$, this shows that
$$
\eqalign{
y-(\deg f)(\deg\L_Y-1)y_0&=[\L_Y^{\deg D}(-f_*D)]\cr
&=\Norm_f\L_X(-D),}
$$
and therefore that the output of the algorithm is indeed $(\Alb
f)(x)$.  Our computational assumptions imply that the running time is
polynomial in~$\deg\L_X$, measured in field operations
in~$k$.\endanalysis

Finally we consider correspondences, i.e.\ diagrams of the form
$$
\correspondence{\;X\;}fg{Y}{Z\rlap,}
$$
where $X$, $Y$ and~$Z$ are proper, smooth, geometrically connected
curves over a field~$k$.  Such a correspondence induces group
homomorphisms
$$
\eqalignno{
\Alb g\circ\Pic f\colon\Pic^0 Y&\to\Pic^0 Z\cr
\noalign{\noindent and}
\Alb f\circ\Pic g\colon\Pic^0 Z&\to\Pic^0 Y.}
$$
Clearly, these can be computed by composing the two algorithms
described above.

\section{Curves over finite fields}

\label{sec:curves-over-finite-fields}

In this section we give algorithms for computing with divisors on a
curve over a finite field.  After some preliminaries, we show how to
compute the Frobenius map on divisors and how to choose uniformly
random divisors of a given degree.  Then we show how to do various
operations in the Picard group of a curve over a finite field, such as
choosing random elements, computing the Frey--R\"uck pairing and
finding a basis of the $l$-torsion for a prime number~$l$.  Many of
the results in this section, especially those
in \S\ref{finding-relations}, \S\ref{kummer-map}
and \S\ref{computing-with-torsion-points}, are variants of work of
Couveignes~\cite{Couveignes: Linearizing torsion classes}.

From now on, we switch from measuring the running time of algorithms
in field operations to measuring it in bit operations.  The usual
field operations in a finite field~$k$ can be done in time polynomial
in~$\log\#k$.

Let $k$ be a finite field of cardinality~$q$, and let $X$ be a
complete, smooth, geometrically connected curve of genus~$g$ over~$k$.
The {\it zeta function\/} of~$X$ is the power series in~$\Z[[t]]$
defined by
$$
\commdiag{
\displaystyle
\llap{$\Zeta_X={}$}
\sum_{D\in\Eff X}t^{\deg D}& \Relbar\joinrel\Relbar&
\displaystyle\sum_{n=0}^\infty(\#\Eff^n X)t^n\cr
\big\Vert & & \big\Vert\cr
\displaystyle\prod_{P\in\PrimeDiv X}{1\over1-t^{\deg P}}&
\Relbar\joinrel\Relbar&
\displaystyle\prod_{d=1}^\infty(1-t^d)^{-\#\PrimeDiv^d X}\rlap.}
$$
Here $\Eff X$ and $\PrimeDiv X$ are the sets of effective divisors and
prime divisors on~$X$, respectively; a superscript denotes the subset
of divisors of the indicated degree.  The following properties of the
zeta function are well known.
\smallskip
\item{(1)} The power series~$\Zeta_X$ can be written as a rational
function
$$
\Zeta_X={L_X\over(1-t)(1-qt)},
\eqnumber{zeta-rational}
$$
where $L_X\in\Z[t]$ is a polynomial of the form
$$
L_X=1+a_1t+\cdots+a_{2g-1}t^{2g-1}+q^gt^{2g}.
$$
\smallskip
\item{(2)} The factorisation of~$L_X$ over the complex numbers has the
form
$$
L_X=\prod_{i=1}^{2g}(1-\alpha_i t),
\eqnumber{zeta-factorisation}
$$
where each $\alpha_i$ has absolute value $\sqrt q$.
\smallskip
\item{(3)} The polynomial~$L_X$ satisfies the functional equation
$$
q^gt^{2g}L_X(1/qt)=L_X(t).
\eqnumber{zeta-functional-equation}
$$
\smallskip

From the definition of~$\Zeta_X$ and from~\eqref{zeta-rational} it is
clear how one can compute the number of effective divisors of a given
degree on~$X$ starting from the polynomial~$L_X$.  We now show how to
extract the number of {\it prime\/} divisors of a given degree
from~$L_X$.  Taking logarithmic derivatives in the definition
of~$\Zeta_X$ and the expression~\eqref{zeta-rational}, we obtain
$$
{\Zeta_X'\over\Zeta_X}={1\over t}\sum_{n=1}^\infty
\Biggl(\sum_{d\mid n}d\cdot\#\PrimeDiv^d X\Biggr)t^n
={L_X'\over L_X}+{1\over1-t}+{q\over1-qt}.
\eqnumber{log-derivative-zeta-function}
$$
Our knowledge of~$L_X$ enables us to compute the coefficients of this
power series.  We can then compute $\#\PrimeDiv^d X$ using the
M\"obius inversion formula.  More explicitly, taking logarithmic
derivatives in the factorisation~\eqref{zeta-factorisation}, we obtain
{\it Newton's identity\/}
$$
L_X'/L_X=-\sum_{n=0}^\infty s_{n+1}t^n,
$$
where the $s_n$ are the power sums
$$
s_n=\sum_{i=1}^{2g}\alpha_i^n\in\Z\quad(n\in\Z).
$$
Expanding the right-hand side of~\eqref{log-derivative-zeta-function}
in a power series and comparing coefficients, we get
$$
\sum_{d\mid n}d\,\#\PrimeDiv^d X=1+q^n-s_n,
$$
or equivalently, by the M\"obius inversion formula,
$$
n\,\#\PrimeDiv^n X=\sum_{d\mid n}\mu(n/d)(1+q^d-s_d),
$$
where $\mu$ is the usual M\"obius function.  Note that this simplifies
to
$$
\#\PrimeDiv^n X=\cases{
1+q-s_1& if $n=1$;\cr
{1\over n}\sum_{d\mid n}\mu(n/d)(q^d-s_d)& if $n\ge2$.}
\eqnumber{formula-number-prime-divisors}
$$

Let $J=\Pic^0_{X/k}$ denote the Jacobian variety of~$X$.  From the
fact that the Brauer group of~$k$ vanishes it follows that the
canonical inclusion
$$
\Pic^0 X\to J(k)
$$
is an equality.  In other words, every rational point of~$J$ can be
identified with a linear equivalence class of $k$-rational divisors of
degree~0.

We note that from the functional
equation~\eqref{zeta-functional-equation} one can deduce that
$$
\#\Eff^n X={q^{1-g+n}-1\over q-1}L_X(1)\quad\hbox{for }n\ge2g,
$$
which in turn is equivalent to ``class number formula''
$$
\#J(k)=\#\Pic^0 X=L_X(1).
\eqnumber{class-number-formula}
$$

\subsection{The Frobenius map}

\label{frobenius-on-divisors}

Let $k$ be a finite field of cardinality~$q$, and let $X$ be a
projective curve over~$k$ in the sense
of~\S\ref{representation-curve}.  We write $d=\deg\L_X$.  Let $\Sym^d
X$ denote the $d$-th symmetric power of~$X$ over~$k$, and let
$\Gr^d\glob(X,\L_X^{\otimes2})$ denote the Grassmann variety of linear
subspaces of codimension~$d$ in the $k$-vector space
$\glob(X,\L_X^{\otimes2})$.  Then we have a commutative diagram
$$
\commdiag{
\Gr^d\glob(X,\L_X^{\otimes2})& \longleftarrow& \Sym^d X\cr
\leftlabel{\Frob_q}\downar& & \downar\rightlabel{\Frob_q}\cr
\Gr^d\glob(X,\L_X^{\otimes2})& \longleftarrow& \Sym^d X\cr}
$$
of varieties over~$k$, where the vertical arrows are the $q$-power
Frobenius morphisms.  Now let $k'$ be a finite extension of~$k$, write
$$
X'=X\times_{\Spec k}\Spec k',
$$
and let $D$ be an effective divisor on~$X'$.  The commutativity of the
above diagram shows that the divisor $\Frob_q(D)$ on~$X'$ can be
computed using the following algorithm.

\algorithm (Frobenius map on divisors). Let $X$ be a projective curve
over a finite field~$k$ of $q$ elements, and let $\Frob_q$ be the
Frobenius map on the set of divisors on~$X$.  Let $k'$ be a finite
extension of~$k$.  Let $X'=X\times_{\Spec k}\Spec k'$, and let
$\L_{X'}$ be the pull-back of the line bundle~$\L_X$ to~$X'$.  Let $i$
be a positive integer, and let $D$ be an effective divisor on~$X'$.
Given the matrix~$M$ of the inclusion map
$$
\glob(X',\L_{X'}^{\otimes i}(-D))\longrightarrow
\glob(X',\L_{X'}^{\otimes i})
$$
with respect to any $k'$-basis of the left-hand side and the
$k'$-basis induced from any $k$-basis of~$\glob(X,\L_X^{\otimes i})$
on the right-hand side, this algorithm outputs the analogous matrix
for the inclusion map
$$
\glob(X',\L_{X'}^{\otimes i}(-\Frob_q(D)))\longrightarrow
\glob(X',\L_{X'}^{\otimes2}).
$$

\label{compute-frobenius-map}

\step Apply the Frobenius automorphism of $k'$ over~$k$ to the
coefficients of the matrix~$M$, and output the result.

\endalgorithm

\analysis It follows from the discussion preceding the algorithm that
the output is indeed equal to $\glob(X',\L_{X'}^{\otimes
i}(-\Frob_q(D)))$.  The algorithm takes $O((\deg\L_X)^2)$ computations
of a $q$-th power of an element in~$k'$.\endanalysis

\subsection{Choosing random prime divisors}

Let $X$ be a projective curve (in the sense
of~\S\ref{representation-curve}) over a finite field.  Our next goal
is to generate random effective divisors of given degree on~$X$.  We
start with an algorithm to generate random prime divisors.  For this
we do not yet need to know the zeta function of~$X$, although we use
its properties in the analysis of the running time of the algorithm.

\algorithm (Random prime divisor). Let $X$ be a projective curve
over a finite field~$k$.  Let $d$ and~$i$ be positive integers such
that
$$
d\le i\deg\L_X-2\genus_X.
$$
Given $d$, $i$ and the $k$-algebra~$S_X^{(2i+2)}$, this algorithm
outputs a uniformly distributed prime divisor~$P$ of degree~$d$
on~$X$, represented as the subspace $\gl(X,\L_X^{\otimes i}(-P))$
of~$\gl(X,\L_X^{\otimes i})$, provided $\PrimeDiv^d X$ is non-empty.
(If $\PrimeDiv^d X=\emptyset$, the algorithm does not terminate.)

\step Choose a non-zero element $s\in\gl(X,\L_X^{\otimes i})$
uniformly randomly, and let $D$ denote the divisor of~$s$.  (In other
words, choose a random hypersurface section of degree~$i$ of~$X$.)

\label{random-prime-divisor-first-step}

\step Compute the set $\Irr^d D$ of (reduced) irreducible components
of~$D$ of degree~$d$ over~$k$ using
Algorithm~\ref{decomposition-divisor}.

\step With probability ${\#\Irr^d D\over\lfloor(i\deg\L_X)/d\rfloor}$,
output a uniformly random element $P\in\Irr^d D$ and stop.

\step Go to step~\ref{random-prime-divisor-first-step}.

\endalgorithm

\label{random-prime-divisor}

\analysis Let $q$ denote the cardinality of~$k$, and let $H$ denote
the set of divisors~$D$ that are divisors of non-zero global sections
of~$\L_X^{\otimes i}$.  By the Riemann--Roch formula, the cardinality
of~$H$ is
$$
\#H={q^{1-g+i\deg\L}-1\over q-1}.
$$
When the algorithm finishes, the probability $p(D,P)$ that a specific
pair $(D,P)$ has been chosen is
$$
\eqalign{
p(D,P)&={1\over\#H}{\#\Irr^d D\over\lfloor(i\deg\L)/d\rfloor}
{1\over\#\Irr D}\cr
&={q-1\over q^{1-g+i\deg\L}-1}{1\over\lfloor(i\deg\L)/d\rfloor}.}
$$
For all prime divisors~$P$ of degree~$d$, the number of $D\in H$ for
which $P$ is in the support of~$D$ is equal to
$$
\#\{D\mid P\in\supp D\}={q^{1-g+i\deg\L-d}-1\over q-1},
$$
so the probability $p(P)$ that a given $P$ is chosen equals
$$
\eqalign{
p(P)&=\#\{D\mid P\in\supp D\}\cdot p(D,P)\cr
&={q^{1-g+i\deg\L-d}-1\over q^{1-g+i\deg\L}-1}
{1\over\lfloor(i\deg\L)/d\rfloor}.}
$$
This is independent of~$P$ and therefore shows that when the algorithm
finishes, the chosen element $P\in\PrimeDiv^d X$ is uniformly
distributed.  Furthermore, the probability~$p$ that the algorithm
finishes in a given iteration is
$$
\eqalign{
p&=\#\PrimeDiv^d X\cdot
{q^{1-g+i\deg\L-d}-1\over q^{1-g+i\deg\L}-1}
{1\over\lfloor(i\deg\L)/d\rfloor}\cr
&={\#\PrimeDiv^d X\over q^d}
{q^{1-g+i\deg\L}-q^d\over q^{1-g+i\deg\L}-1}
{1\over\lfloor(i\deg\L)/d\rfloor}\cr
&\ge{\#\PrimeDiv^d X\over q^d}(1-q^{-1-\genus_X}){d\over i\deg\L}.}
$$
We claim that the expected running time is polynomial in $\deg\L$, $i$
and~$\log q$, under the assumption that $\#\PrimeDiv^d X\ne\emptyset$.
we distinguish two cases:
$$
q^{d/2}<2\sigma^0(d)(2\genus_X+1)\quad\hbox{and}\quad
q^{d/2}\ge2\sigma^0(d)(2\genus_X+1),
$$
where $\sigma^0(d)$ denotes the number of positive divisors of~$d$.
In the first case, we see that
$$
p>(2\sigma^0(d)(2\genus_X+1))^2(1-q^{-1-\genus_X}){d\over i\deg\L},
$$
which shows that $1/p$ is bounded by a polynomial in $\deg\L$ and~$i$,
In the second case, we deduce
from~\eqref{formula-number-prime-divisors} the following estimate
for~$\#\PrimeDiv^d X$:
$$
\eqalign{
|d\#\PrimeDiv^d X-q^d|&\le\sum_{\textstyle{e\mid d\atop e\ne d}}q^e
+\sum_{e\mid d}|s_e|\cr
&\le(\sigma^0(d)-1)q^{d/2}+\sigma^0(d)\cdot 2\genus_X q^{d/2}\cr
&<\sigma^0(2)(2\genus_X+1)q^{d/2}\cr
&\le{1\over2}q^d},
$$
so that $\#\PrimeDiv^d X>q^d/(2d)$, and hence
$$
p>{1-q^{-1-\genus_X}\over2i\deg\L}.
$$
In both cases we conclude that the expected running time is bounded by
a polynomial in $\deg\L$, $i$ and~$\log q$.\endanalysis

\subsection{Choosing random divisors}

As before, let $X$ be a projective curve over a finite field~$k$.
From now on we assume that we know the zeta function of~$X$, or
equivalently the polynomial~$L_X$.

Below we will give an algorithm for generating uniformly random
effective divisors of a given degree on the curve~$X$.  These divisors
will be built up from prime divisors, so it will be useful to speak of
the {\it decomposition type\/} of an effective divisor~$D$.  This is
the sequence of integers $(l_1,l_2,\ldots)$, where $l_d$ is the number
of prime divisors of degree~$d$ (counted with multiplicities)
occurring in~$D$.

One of the ingredients is the concept of {\it $m$-smooth\/} divisors
and decomposition types.  An $m$-smooth divisor is a linear
combination of prime divisors whose degrees are at most~$m$, and an
$m$-smooth decomposition type of degree~$n$ is an
$m$-tuple~$(l_1,\ldots,l_m)$ such that $\sum_{d=1}^m l_dd=n$.  For
every $m$-smooth effective divisor~$D$ of degree~$n$, we may view the
decomposition type of~$D$ as an $m$-smooth decomposition type, since
only its first $m$ coeffients are non-zero.

The algorithm that we will describe takes as input the degree~$n$ as
well as a positive integer~$m$, and outputs a uniformly random
$m$-smooth effective divisor of degree~$n$.  Clearly, all effective
divisors of degree~$n$ are $n$-smooth, so that the algorithm can be
used with $m=n$ to produce uniformly random effective divisors of
degree~$n$.

The first step is to generate the decomposition type of a uniformly
random $m$-smooth effective divisor of degree~$n$.  The method we use
for doing this is described by Diem in~\citex{Diem}{page~150} and
in~\citex{Diem: article}.  Diem's algorithm works by recursion on~$m$.

For every $m\ge1$, we write $\Eff^n_{\le m}X$ for the set of
$m$-smooth effective divisors~$D$ of degree~$n$.  Furthermore, for
$l\ge0$ and $m\ge1$ we write $\Eff^{lm}_{=m}X$ for the set of divisors
of degree~$lm$ that are linear combinations of prime divisors of
degree~$m$.  We note that the set~$\Eff^n_{\le m}X$ can be decomposed
as
$$
\Eff^n_{\le m}X=\cases{\displaystyle\strut\Eff^n_{=1}X& if $m=1$;\cr
\displaystyle\bigsqcup_{l=0}^{\lfloor n/m\rfloor}
\Eff^{lm}_{=m}X\times\Eff^{n-lm}_{\le m-1}X& if $m\ge2$.}
\eqnumber{decomposition-smooth-divisors}
$$
The cardinality of~$\Eff^{lm}_{=m}X$ equals the number of ways to
choose $l$ elements from the set~$\PrimeDiv^m X$ with repeats.  For
this we have the well-known formula
$$
\#\Eff^{lm}_{=m}X={\#\PrimeDiv^m X-1+l\choose l}.
\eqnumber{recursion-effective-divisors-1}
$$
Furthermore, from the
description~\eqref{decomposition-smooth-divisors} of~$\Eff^n_{\le m}X$
we see that
$$
\#\Eff^n_{\le m}X=\cases{
\displaystyle\strut\#\Eff^n_{=1}X& if $m=1$;\cr
\displaystyle\sum_{l=0}^{\lfloor n/m\rfloor}
\#\Eff^{lm}_{=m}X\cdot\#\Eff^{n-lm}_{\le m-1}X&
if $m\ge2$.}
\eqnumber{recursion-effective-divisors-2}
$$
From these relations we can compute $\#\Eff^n_{\le m}X$ recursively,
starting from the numbers~$\#\PrimeDiv^d X$ for~$1\le d\le m$.  An
alternative way to describe these recurrence relations is to use
generating functions; see Diem~\citex{Diem}{page~149} or \citex{Diem:
article}{Lemma~3.14}.

In order to generate decomposition types of uniformly random
$m$-smooth divisors of degree~$n$, we define a probability
distribution~$\mu^n_m$ on the set of $m$-smooth decomposition types of
degree~$n$ by defining $\mu^n_m(l_1,\ldots,l_m)$ as the probability
that a uniformly randomly chosen effective $m$-smooth divisor of
degree~$n$ has decomposition type~$(l_1,\ldots,l_m)$.  The algorithm
now works as follows.  We first select an
integer~$l_m\in\{0,1,\ldots,\lfloor n/m\rfloor\}$---the number of
prime divisors of degree~$m$ (counted with multiplicities) occurring
in the decomposition---according to the marginal
distribution~$\nu^n_m$ of the $m$-th coordinate.  We then apply the
algorithm recursively with $(n-l_mm,m-1)$ in place of~$(n,m)$.

The marginal distribution~$\nu^n_m$ of the coordinate~$l_m$ in a
$m$-tuple~$(l_1,\ldots,l_m)$ distributed according to~$\mu^n_m$ is the
following.  If $m=1$, then $l_1=n$ with probability~1.  When $m\ge2$,
the probability that $l_m$ equals a given $l\in\{0,1,\dots,\lfloor
n/m\rfloor\}$ is
$$
\nu^n_m(l)={\#\Eff^{lm}_{=m}X\cdot\#\Eff^{n-lm}_{\le m-1}X\over
\#\Eff^n_{\le m}X}\quad(0\le l\le\lfloor n/m\rfloor).
\eqnumber{distribution-l}
$$
Once we have computed $\#\Eff^n_{\le m}X$, as well as
$\#\Eff^{lm}_{=m}$ and $\#\Eff^{n-lm}_{\le m-1}X$ for $0\le
l\le\lfloor n/m\rfloor$ (using \eqref{formula-number-prime-divisors},
\eqref{recursion-effective-divisors-1}
and~\eqref{recursion-effective-divisors-2}), it is straightforward to
generate a random~$l_m\in\{0,1,\ldots,\lfloor n/m\rfloor\}$
distributed according to~$\nu^n_m$.  Namely, we subdivide the interval
$$
I=\{0,1,\ldots,\#\Eff^n_{\le m}X-1\}
$$
into $\lfloor n/m\rfloor+1$ intervals~$I_l$, with $0\le l\le\lfloor
n/m\rfloor$ and each $I_l$ having length
$\#\Eff^{lm}_{=m}X\cdot\#\Eff^{n-lm}_{\le m-1}X$, we generate a
uniformly random element $x\in I$, and we select the unique $l$ such
that $x\in I_l$.

\algorithm (Decomposition type of a random divisor).  Given the
polynomial~$L_X$ for a curve~$X$ over a finite field and integers
$n\ge0$ and $m\ge1$, this algorithm outputs a random $m$-smooth
decomposition type~$(l_1,\ldots,l_m)$ of degree~$n$, distributed
according to the distribution~$\mu^n_m$.

\label{random-decomposition-type}

\step If $m=1$, output the 1-tuple $(n)$ and stop.

\step Choose a random element $l_m\in\{0,1,\ldots,\lfloor
n/m\rfloor\}$ according to the distribution~$\nu^n_m$
from~\eqref{distribution-l}.

\step Call the algorithm recursively with $(n-l_m m,m-1)$ in place
of~$(n,m)$ to obtain an $(m-1)$-smooth decomposition
type~$(l_1,\ldots,l_{m-1})$ of degree~$n-l_m m$.

\step Output the $m$-tuple $(l_1,\ldots,l_m)$.

\endalgorithm

\analysis The correctness of the algorithm follows from the above
discussion.  It is straightforward to check that it runs in time
polynomial in $\genus_X$, $\log\#k$, $n$ and~$m$.\endanalysis

The preceding algorithm reduces our problem to generating random
linear combinations of $l$ prime divisors of a given degree~$d$.  In
other words, we have to pick a random {\it multiset\/} of
cardinality~$l$ from~$\PrimeDiv^d X$.  This can be done using the
following algoritm.  I thank Claus Diem for pointing out this method
to me, which is much simpler than the one I had in mind originally.

\algorithm (Random multiset). Let $S$ be a finite non-empty set of
known cardinality.  Suppose we have algorithms to pick uniformly
random elements of~$S$ and to decide whether two such elements are
equal.  Given a non-negative integer~$l$, this algorithm outputs a
uniformly random multiset of $l$ elements from~$S$.

\step Generate a uniformly random subset $\{x_1,\ldots,x_l\}$ of
$\{1,2,\ldots,l+\#S-1\}$, with $x_1<x_2<\ldots<x_l$.

\step Define a multiset $(y_1,\ldots,y_l)$ of $l$ elements
from~$\{0,1,\ldots,\#S-1\}$ by $y_i=x_i-i$; then $y_1\le
y_2\le\ldots\le y_l$.

\step For each $i$ with $1\le i\le l$, let $a_i$ be the number of
elements of~$\{0,1,\ldots,\#S-1\}$ that occur with multiplicity~$i$
in~$(y_1,\ldots,y_l)$.

\step Generate a uniformly random sequence
$$
\displaylines{
s_1^1,s_2^1,\ldots,s_{a_1}^1,\cr
s_1^2,s_2^2,\ldots,s_{a_2}^2,\cr
\vdots\cr
s_1^l,s_2^l,\ldots,s_{a_l}^l}
$$
of $a_1+a_2+\cdots+a_l$ distinct elements of~$S$.

\step Output the multiset consisting of the elements $s_i^j$ of~$S$,
where $s_i^j$ occurs with multiplicity~$j$.

\endalgorithm

\label{random-divisor-finite-set}

\analysis By construction, the multiset $(y_1,\ldots,y_l)$ of $l$
elements from~$\{0,1,\ldots,\#S-1\}$ is uniformly random, so the
``multiplicity vector'' $(a_1,\ldots,a_l)$ is the same as that of a
uniformly random multiset of $l$ elements from~$S$.  The multiset
generated in the last step is uniformly random among the multisets
with this ``multiplicity vector''.  This implies that the result is a
uniformly random multiset of $l$ elements from~$S$, as
required.\endanalysis

Combining Algorithms \ref{random-prime-divisor},
\ref{random-decomposition-type} and~\ref{random-divisor-finite-set},
we obtain the following algorithm to generate a uniformly random
effective divisor of a given degree.

\algorithm (Random divisor). Let $X$ be a projective curve over a
finite field~$k$.  Given positive integers $m$ and~$i$, an integer~$n$
satisfying
$$
0\le n\le i\deg\L_X-2\genus_X,
$$
the graded $k$-algebra~$S_X^{(2i+2)}$ and the polynomial~$L_X$, this
algorithm outputs a uniformly random $m$-smooth effective divisor~$D$
of degree~$n$ on~$X$, represented as the
subspace~$\gl(X,\L_X^{\otimes i}(-D))$ of~$\gl(X,\L_X^{\otimes
i})$.

\label{random-divisor-curve}

\step Generate a random $m$-smooth decomposition
type~$(l_1,\ldots,l_m)$ of degree~$n$ using
Algorithm~\ref{random-decomposition-type}.

\step For $d=1,\ldots,m$, generate a uniformly random linear
combination $D_d$ of $l_d$ prime divisors of degree~$d$ on~$X$ using
Algorithm~\ref{random-divisor-finite-set} (with $S=\PrimeDiv^d X$, and
$l=l_d$), where we use Algorithm~\ref{random-prime-divisor} to
generate random elements of~$\PrimeDiv^d X$.

\step Compute the subspace $\gl(X,\L_X(-D))$ for the divisor
$D=D_1+\cdots+D_m$ using the addition algorithm described
in~\S\ref{representation-divisors}, and output $\gl(X,\L_X(-D))$.

\endalgorithm

\analysis It follows from the above discussion that the algorithm
outputs a uniformly random $m$-smooth divisor of degree~$n$ on~$X$.
The running time is clearly polynomial in $m$, $n$, $i$ and~$\deg\L_X$
(measured in field operations in~$k$).\endanalysis

\remark In practice, the following method for picking a random
effective divisor of degree~$n$ is faster, but does not give a
uniformly distributed output.  We first choose a uniformly random
non-zero section~$s$ of~$\glob(X,\L^{\otimes i})$, where $i$ is a
non-negative integer such that
$$
i\deg\L-n\ge 2g+1.
$$
Then if the set of effective divisors~$D$ of degree~$n$ with $D\le\div
s$ is non-empty, we pick a uniformly random element from it; otherwise
we keep going with a different section~$s$.

\subsection{The Frobenius endomorphism of the Jacobian}

\label{frobenius-endomorphism}

As before, let $k$ be a finite field of cardinality~$q$, and let $X$
be a proper, smooth and geometrically connected curve over~$k$.  Let
$J$ be the Jacobian variety of~$X$, and let $\Frob_q$ denote the
Frobenius endomorphism of~$J$; is an isogeny of degree~$q$.  The
Rosati dual of~$\Frob_q$ is called the {\it Verschiebung\/} and
denoted by~$\Ver_q$.  The Albanese and Picard maps associated to the
Frobenius morphism on~$X$ are the endomorphisms $\Frob_q$ and~$\Ver_q$
of~$J$, respectively.

Then we have a commutative diagram
$$
\commdiag{
\Sym^d X& \longrightarrow& J\cr
\leftlabel{\Frob_q}\downar& & \downar\rightlabel{\Frob_q}\cr
\Sym^d X& \longrightarrow& J\cr}
$$
of varieties over~$k$, where the vertical arrows are the $q$-power
Frobenius morphisms.  This shows that the Frobenius endomorphism
of~$J$ is equal to the endomorphism $\Alb(\Frob_q)$ induced by the
Frobenius map on~$X$ via Albanese functoriality.

Write $X'=X\times_{\Spec k}\Spec k'$.  The results
of~\S\ref{frobenius-on-divisors} now imply that for any finite
extension $k'$ of~$k$, the endomorphism $\Frob_q$
of~$J(k')=\Pic^0(X')$ can be computed by applying
Algorithm~\ref{compute-frobenius-map} to any subspace
$\glob(X',\L_{X'}^{\otimes2}(-D))$ of the $k'$-vector space
$$
\glob(X',\L_{X'}^{\otimes 2})\cong
k'\otimes_k\glob(X,\L_X^{\otimes 2})
$$
where $D$ is an effective divisor of degree~$\deg\L_X$ on~$X'$ such
that $\L_{X'}(-D)$ represents $x$.

If $O$ is a $k$-rational point of~$X$, then we can compute the trace
map
$$
\tr_{k'/k}\colon\Pic^0 X'\to \Pic^0 X
$$
in the following way.  For $x\in\Pic^0 X'$, we compute a subspace
of~$\glob(X',\L_{X'}^{\otimes2})$ representing the element
$$
y=\sum_{i=0}^{[k':k]}\Frob_qx\in\Pic^0 X'.
$$
Now $y$ is in fact the image of the element~$\tr_{k'/k}x\in\Pic^0 X$
under the inclusion $\Pic^0 X\to\Pic^0 X'$, so we can apply
Algorithm~\ref{descent} to find a subspace
of~$\glob(X,\L_X^{\otimes2})$ representing $\tr_{k'/k}x$.

In \S\ref{picard-albanese}, the problem of computing the Albanese map
for a finite morphism of curves was reduced to the problem of compute
trace maps.  Since we can solve the latter problem, we can also solve
the former.

\subsection{Picking random elements of the Picard group}

\label{random-elements-picard-group}

The next problem we will study is that of picking uniformly random
elements in the finite Abelian group~$J(k)=\Pic^0 X$.  We recall
from~\S\ref{computing-in-picard-group} that in the medium model of the
Picard group, the class of a line bundle~${\cal M}$ of degree~0 is
represented by an effective divisor~$D$ of degree~$\deg\L$ such that
${\cal M}\cong\L(-D)$.  Consider the map
$$
\eqalign{
\Eff^{\deg\L}X&\to\Pic^0 X\cr
D&\mapsto[\L(-D)].}
$$
It follows from the Riemann--Roch theorem and the fact that
$\deg\L\ge2\genus_X-1$ that all fibres of this map have cardinality
${q^{1-g+\deg\L}-1\over q-1}$.  This means that to pick a uniformly
random element of~$\Pic^0 X$ it suffices to pick a uniformly random
divisor of degree~$\deg\L$.  A method for doing this is given by
Algorithm~\ref{random-divisor-curve}, provided that we know
$S_X^{(6)}$.

\subsection{Computing Frey--R\"uck pairings}

\label{computing-frey-rueck-pairings}

Let $n$ be a positive integer.  We assume $k$ contains a primitive
$n$-th root of unity; this is equivalent to
$$
n\mid\#k^\times=q-1
$$
and implies that $n$ is not divisible by the characteristic of~$k$.

Let $X$ be a complete, smooth, geometrically connected curve over~$k$,
and let $J$ be its Jacobian variety.  The {\it Frey--R\"uck pairing\/}
of order~$n$ on~$J(k)=\Pic^0 X$, often also referred to as the {\it
Tate--Lichtenbaum pairing\/}, is the bilinear map
$$
[\blank\mathord,\blank]_n\colon
J[n](k)\times J(k)/nJ(k)\to\mu_n(k)
$$
defined as follows (see Frey and R\"uck~\cite{Frey-Rueck} or
Schaefer~\cite{Schaefer}).  Let $x$ and~$y$ be elements of~$J(k)$ such
that $nx=0$.  Choose divisors $D$ and~$E$ such that $x$ and~$y$ are
represented by the line bundles $\O_X(D)$ and~$O_X(E)$, respectively,
and such that the supports of $D$ and~$E$ are disjoint.  By
assumption, there exists a rational function~$f$ on~$X$ such that
$nD=\div(f)$; now $[x,y]_n$ is defined as
$$
[x,y]_n=f(E)^{\#k^\times/n}.
$$
Here $f(E)$ is defined on $\bar k$-valued points (where $\bar k$ is an
algebraic closure of~$k$) by function evaluation, and then extended to
the group of divisors on~$X_{\bar k}$, by linearity in the sense that
$$
f(E+E')=f(E)\cdot f(E').
$$
It is known that the Frey--R\"uck pairing is {\it perfect\/} in the
sense that it induces isomorphisms
$$
\eqalignno{
J[n](k)&\isom\Hom(J(k)/nJ(k),\mu_n(k))\cr
\noalign{\noindent and}
J(k)/nJ(k)&\isom\Hom(J[n](k),\mu_n(k))
}
$$
of Abelian groups.

Let us now give a slightly different interpretation of~$f(E)$ that
brings us in the right situation to compute $[x,y]_n$.  We consider an
arbitrary non-zero rational function~$f$ and an arbitrary divisor~$E$
such that the divisors
$$
D=\div(f)
$$
and~$E$ have disjoint supports.  Since $f(E)$ is by definition linear
in~$E$, it suffices to consider the case where $E$ is an effective
divisor.  As in~\S\ref{norm-functor-effective-divisors}, we write
$$
j_E\colon E\to X
$$
for the closed immersion of~$E$ into~$X$, and if ${\cal M}$ is a line
bundle on~$X$ we abbreviate
$$
\Norm_{E/k}{\cal M}=\Norm_{E/k}(j_E^*{\cal M}).
$$
Since $D$ and~$E$ have disjoint supports, we have a canonical
trivialisation
$$
t_D\colon k\cong\Norm_{E/k}\O_X\isom\Norm_{E/k}\O_X(D).
$$
On the other hand, multiplication by~$f$ induces an isomorphism
$$
\Norm_{E/k}f\colon\Norm_{E/k}\O_X(D)\isom\Norm_{E/k}\O_X\cong k.
$$
of one-dimensional $k$-vector spaces.  We claim that the composed
isomorphism
$$
k\isomorphism{t_D}\Norm_{E/k}\O_X(D)\isomorphism{\Norm_{E/k}f}k
\eqnumber{evaluate-function-on-divisor}
$$
is multiplication by~$f(E)$.  This is true in the case where $E$ is a
single point, since then $\Norm_{E/k}$ is (canonically isomorphic to)
the identity functor.  We deduce the general case from this by
extending the base field to an algebraic closure of~$k$ and using the
fact that both $f(E)$ and the norm functor are linear in~$E$.  For the
latter claim, we refer to Deligne~\citex{SGA 4 III}{expos\'e~XVII, \no
6.3.27}.

\remark The isomorphism~\eqref{evaluate-function-on-divisor} could be
taken as a {\it definition\/} of~$f(E)$ for effective divisors~$E$.

\proclaim Lemma.  Let $x$ and~$y$ be elements of~$J(k)$ with $nx=0$,
let ${\cal M}$ be a line bundle representing $x$, and let $E^+$
and~$E^-$ be effective divisors such that $\O_X(E^+-E^-)$ represents
$y$.  (In particular, ${\cal M}$ has degree~0, and $E^+$ and~$E^-$
have the same degree.)  For any pair of trivialisations
$$
t^\pm\colon k\isom\Norm_{E^\pm/k}{\cal M}
$$
of $k$-vector spaces and any trivialisation
$$
s\colon\O_X\isom{\cal M}^{\otimes n}
$$
of line bundles on~$X$, the isomorphism
$$
k\isomorphism{(t^+)^n}\Norm_{E^+/k}{\cal M}^{\otimes n}
\isomorphism{\Norm_{E^+/k}s^{-1}}
k\isomorphism{\Norm_{E^-/k}s}\Norm_{E^-/k}{\cal M}^{\otimes n}
\isomorphism{(t^-)^{-n}}k
$$
is multiplication by an element of~$k^\times$ whose
$(\#k^\times/n)$-th power equals $[x,y]_n$.

\label{frey-rueck-via-norm}

\noindent
(We have implicitly used the isomorphisms $\Norm_{E^\pm/k}({\cal
M}^{\otimes n})\cong(\Norm_{E^\pm/k}{\cal M})^{\otimes n}$ expressing
the linearity of~$\Norm_{E/k}$, and denoted both sides of the
isomorphism by $\Norm_{E^\pm/k}{\cal M}^{\otimes n}$.)

\medskip\nobreak

\proof We fix a non-zero rational section~$h$ such that the divisor
$$
D=\div h
$$
is disjoint with~$E^\pm$.  Then we have canonical trivialisations
$$
t_D^\pm\colon k\isom\Norm_{E^\pm/k}\O_X(D)
$$
as above.  Composing these with the isomorphism
$$
\Norm_{E^\pm/k}h\colon\Norm_{E^\pm/k}\O_X(D)
\isom\Norm_{E^\pm/k}{\cal M}
$$
induced by multiplication by~$h$ gives trivialisations
$$
t_h^\pm=\Norm_{E^\pm/k}h\circ t_D\colon k\isom\Norm_{E^\pm/k}{\cal M}.
$$
Now consider any isomorphism
$$
s\colon\O_X\isom{\cal M}^{\otimes n}
$$
of line bundles on~$X$, and define
$$
f=s^{-1}\circ h^n\colon\O_X(nD)\isom\O_X;
$$
then $f$ can be viewed as a rational function with divisor~$nD$.  We
now have commutative diagrams
$$
\commdiag{
k& \isomorphism{(t_D^\pm)^n}& \Norm_{E^\pm/k}\O_X(nD)&
\isomorphism{\Norm_{E^\pm/k}f}& k\cr
\bigm\|& & \leftlabel\sim\downar\rightlabel{\Norm_{E^\pm/k}h^n}& & \bigm\|\cr
k& \isomorphism{(t_h^\pm)^n}& \Norm_{E^\pm/k}{\cal M}^{\otimes n}&
\isomorphism{\Norm_{E^\pm/k}s^{-1}}& k\rlap.}
$$
As we saw above, the top row is multiplication by~$f(E^\pm)$; by the
commutativity of the diagram, the same holds for the bottom row.
Finally, we note that replacing $t_h^\pm$ by {\it any\/} pair of
trivialisations
$$
t^\pm\colon k\isom\Norm_{E^\pm/k}{\cal M}
$$
changes the isomorphism in the bottom row of the above diagram by some
$n$-th power in~$k^\times$.  This implies that the isomorphism
$$
k\isomorphism{(t^\pm)^n}\Norm_{E^\pm/k}{\cal M}^{\otimes n}
\isomorphism{\Norm_{E^\pm/k}s^{-1}}k
$$
equals multiplication by an element of~$k^\times$ whose $(\#k/n)$-th
power is $f(E^\pm)^{\#k^\times/n}$.  The lemma follows from this by
the definition of~$[x,y]_n$.\endproof

Lemma~\ref{frey-rueck-via-norm} reduces the problem of computing the
Frey--R\"uck pairing of order~$n$ to the following: given a line
bundle~${\cal M}$ such that ${\cal M}^{\otimes n}$ is trivial, find an
isomorphism
$$
s\colon\O_X\isom{\cal M}^n,
$$
and, given moreover an effective divisor~$E$ and a trivialisation
$$
t\colon k\isom\Norm_{E/k}{\cal M},
$$
compute the isomorphism
$$
I^E_{s,t}\colon k\isomorphism{t^n}\Norm_{E/k}{\cal M}^{\otimes n}
\isomorphism{\Norm_{E/k}s^{-1}}k.
\eqnumber{isom-Est}
$$

We assume that the curve~$X$ is given by a projective embedding via a
line bundle~$\L$ as in~\S\ref{representation-curve}.  We will describe
an algorithm to compute isomorphisms of the type~$I^E_{s,t}$, based on
Khuri-Makdisi's algorithms for computing with divisors on~$X$.
Suppose we are given a line bundle~${\cal M}$ of degree~0 such that
${\cal M}^{\otimes n}$ is trivial and an effective divisor~$E$.  For
simplicity, we assume that $\deg E=\deg\L$.  As
in~\S\ref{representation-divisors}, we represent the class of~${\cal
M}$ in~$J(k)$ by the subspace~$\glob(X,\L^{\otimes2}(-D))$
of~$\glob(X,\L^{\otimes2})$, where $D$ is any effective divisor of
degree~$\deg\L$ (not necessarily disjoint from~$E$) such that
$$
{\cal M}\cong\L(-D).
$$
Likewise, we represent $E$ as the
subspace~$\glob(X,\L^{\otimes2}(-E))$ of~$\glob(X,\L^{\otimes2})$.

First, we will describe a construction of a trivialisation
$$
s\colon\O_X\isom\L(-D)^{\otimes n}.
$$
For this we fix an anti-addition chain $(a_0,a_1,\ldots,a_m)$ for~$n$,
as described in~\S\ref{computing-in-picard-group}.  In particular, for
each $l$ with $2\le l\le m$ we are given $i(l)$ and $j(l)$ in
$\{0,1,\ldots,l-1\}$ such that
$$
a_l=-a_{i(l)}-a_{j(l)}.
$$
We fix any non-zero global section~$u$ of~$\L$, and we put
$$
D_0=\div(u),\quad D_1=D.
$$
For $l=2$, 3, \dots, $m$, we iteratively apply Algorithm~\ref{addflip}
to $D_{i(l)}$ and~$D_{j(l)}$; this gives an effective divisor~$D_l$ of
degree~$\deg\L$ and a global section~$s_l$ of~$\L^{\otimes3}$ such
that the line bundle $\L^{\otimes3}(-D_l-D_{i(l)}-D_{j(l)})$ is
trivial and
$$
\div(s_l)=D_l+D_{i(l)}+D_{j(l)}.
$$
We recursively define rational sections $h_1$, $h_2$, \dots, $h_m$
of~$\L^{\otimes(a_l-1)}$ by
$$
h_l=\cases{
u^{-1}& for $l=0$;\cr
1& for $l=1$;\cr
(h_{i(l)}h_{j(l)}s_l)^{-1}& for $l=2$, 3, \dots, $m$.}
$$
Then it follows immediately that each $h_l$ has divisor $a_l D-D_l$.
In particular, since $\L(-D)^{\otimes n}$ is trivial, so is $\L(-D_m)$
and Algorithm~\ref{zero-test} provides us with a global section~$v$
of~$\L$ such that
$$
\div(v)=D_m.
$$
The rational section
$$
s=h_m v
$$
of~$\L^{\otimes n}$ has divisor~$nD$ and hence induces an isomorphism
$$
s\colon\O_X\isom\L(-D)^{\otimes n}.
$$

Next, we assume that an effective divisor~$E$ has been given.  We
assume for simplicity that $\deg E=\deg\L$.  We fix bases of the
following $k$-vector spaces:
$$
\displaylines{
\glob(E,\L^{\otimes2});\cr
\glob(E,\L^{\otimes 3}(-D_l))\hbox{ for }1\le l\le m;\cr
\glob(E,\L^{\otimes 4}(-D_{i(l)}-D_{j(l)}))\hbox{ for }2\le l\le m.}
$$
In addition, we fix a $k$-basis of~$\glob(E,\L^{\otimes3}(-D_0))$ by
defining it as the image of the chosen basis of~$\glob(E,\L^{\otimes
2})$ under the multiplication map
$$
u\colon\glob(E,\L^{\otimes2})\isom\glob(E,\L^{\otimes3}(-D_0)).
$$
For $0\le l\le m$ we define a trivialisation
$$
\eqalign{
t_l\colon k&\isom\Norm_{E/k}\L(-D_l)\cr
&\isom\Hom_k\bigl(\det_k\glob(E,\L^{\otimes2}),
\det_k\glob(E,\L^{\otimes3}(-D_l))\bigr)}
$$
using the given bases of~$\glob(E,\L^{\otimes2})$ and
$\glob(E,\L^{\otimes3}(-D_l))$, and we define an element $\gamma_l$
of~$k^\times$ by requiring that the diagram
$$
\commdiag{
k& \isomorphism{t_l}& \Norm_{E/k}\L(-D_l)\cr
\leftlabel{\gamma_l}\downar\rightlabel\sim& &
\leftlabel\sim\downar\rightlabel{h_l}\cr
k& \isomorphism{t^{a_l}}& \Norm_{E/k}\L(-D)^{\otimes a_l}}
$$
be commutative.  For $2\le l\le m$, we define a trivalisation
$$
t_l'\colon k\isom\Norm_{E/k}\L^{\otimes2}(-D_{i(l)}-D_{j(l)})
$$
by \eqref{norm-det-pm} using the given bases of
$\glob(E,\L^{\otimes2})$ and
$\glob(E,\L^{\otimes4}(-D_{i(l)}-D_{j(l)}))$, and a trivialisation
$$
t_l''\colon k\isom\Norm_{E/k}\L^{\otimes3}(-D_l-D_{i(l)}-D_{j(l)})
$$
by \eqref{norm-det-pm} using the given bases of
$\glob(E,\L^{\otimes2})$ and
$\glob(E,\L^{\otimes5}(-D_l-D_{i(l)}-D_{j(l)}))$.

\algorithm (Compute isomorphisms of the form $I^E_{s,t}$).  Let $X$ be
a projective curve over a field~$k$, let $D$ and~$E$ be effective
divisors of degree~$\deg\L$ on~$X$, and let $n$ be a positive integer
such that $\L(-D)^{\otimes n}$ is trivial.  Given the $k$-algebra
$S_X^{(7)}$, an anti-addition chain $(a_0,a_1,\ldots,a_m)$ for~$n$, a
global section~$u$ of~$\L$, effective divisors $D_0$, $D_1$, \dots,
$D_m$, global sections $s_2$, \dots, $s_m$ of~$\L^3$ such that
$$
D_0=\div(u),D_1=D\quad\hbox{and}\quad
\div(s_l)=D_l+D_{i(l)}+D_{j(l)}\hbox{ for }2\le l\le m
$$
and a global section $v$ of the trivial line bundle~$\L(-D_m)$, this
algorithm outputs the isomorphism $I^E_{s,t}$ defined
by~\eqref{isom-Est}, where $s$ is defined using the given data, and
where $t$ is chosen by the algorithm.  (This means that the output of
the algorithm is an element of~$k^\times$ defined up to $n$-th powers
in~$k^\times$.)

\step Put $\gamma_0=\gamma_1=1$.

\step For $l=2$, 3, \dots, $m$:

\plusindent

\step Using Algorithm~\ref{algorithm-linearity-norm}, compute the
elements $\lambda_l^{(1)}$ and~$\lambda_l^{(2)}$ of~$k^\times$ such
that the diagrams
$$
\commdiag{
k& \isomorphism{t_{i(l)}\otimes t_{j(l)}}&
\Norm_{E/k}\L(-D_{i(l)})\otimes\Norm_{E/k}\L(-D_{j(l)})\cr
\leftlabel{\lambda_l^{(1)}}\downar\rightlabel\sim& &
\downar\rightlabel\sim\cr
k& \isomorphism{t_l'}& \Norm_{E/k}\L^{\otimes2}(-D_{i(l)}-D_{j(l)})}
$$
and
$$
\commdiag{
k& \isomorphism{t_l\otimes t_l'}& \Norm_{E/k}\L(-D_l)
\otimes\Norm_{E/k}\L^{\otimes2}(-D_{i(l)}-D_{j(l)})\cr
\leftlabel{\lambda_l^{(2)}}\downar\rightlabel\sim& &
\downar\rightlabel\sim\cr
k& \isomorphism{t_l''}&
\Norm_{E/k}\L^{\otimes3}(-D_l-D_{i(l)}-D_{j(l)})}
$$
are commutative.  Define $\lambda_l=\lambda_l^{(1)}\lambda_l^{(2)}$.

\step Compute $\sigma_l\in k^\times$ as the determinant of the matrix
of the isomorphism
$$
s_l\colon\glob(E,\L^{\otimes2})\isom
\glob(E,\L^{\otimes5}(-D_l-D_{i(l)}-D_{j(l)}))
$$
with respect to the given bases.

\step Put $\displaystyle\gamma_l={\lambda_l\over
\sigma_l\gamma_{i(l)}\gamma_{j(l)}}$.

\minusindent

\step Compute $\delta\in k^\times$ as the determinant of the matrix of
the isomorphism
$$
v\colon\glob(E,\L^2)\isom\glob(E,\L^3(-D_m))
$$
with respect to the given bases.

\step Output the element $\displaystyle{1\over\gamma_m\delta}\in
k^\times$.

\endalgorithm

\label{algorithm-Est}

\analysis The definitions of $\lambda_l$ and~$\sigma_l$ given in the
algorithm imply that the diagram
$$
\commdiag{
k& \isomorphism{t_l\otimes t_{i(l)}\otimes t_{j(l)}}&
\Norm_{E/k}\L(-D_l)\otimes\Norm_{E/k}\L(-D_{i(l)})\otimes
\Norm_{E/k}\L(-D_{j(l)})\cr
\leftlabel{\lambda_l}\downar\rightlabel\sim& &
\downar\rightlabel\sim\cr
k& \isomorphism{t_l''}&
\Norm_{E/k}\L^{\otimes3}(-D_l-D_{i(l)}-D_{j(l)})}
$$
is commutative and that the isomorphism
$$
k\isomorphism{s_l}\Norm_{E/k}\L^{\otimes3}(-D_l-D_{i(l)}-D_{j(l)})
\isomorphism{(t_l'')^{-1}} k
$$
is multiplication by~$\sigma_l$.

The recursive definition of the $h_l$ implies that the recurrence
relation between the $\gamma_l$ is as stated in the algorithm.
Namely, it follows from the definition of~$D_0$, from the special
choice of basis of~$\glob(E,\L^{\otimes3}(-D_0))$ and from the fact
that $t_1=t$ that
$$
\gamma_0=\gamma_1=1.
$$
Furthermore, the definitions of $h_l$, $\gamma_l$, $\gamma_{i(l)}$,
$\gamma_{j(l)}$ and the properties of $\lambda_l$ and~$\sigma_l$ that
we have just proved imply that
$$
\gamma_l={\lambda_l\over\sigma_l\gamma_{i(l)}\gamma_{j(l)}}
\quad\hbox{for }l=2,3,\ldots,m.
$$
Finally, it follows from the definitions of $s$, $\gamma_m$ and the
isomorphism~$I^E_{s,t}$ from~\eqref{isom-Est} that the relation
between $v$, $t_m$, $\gamma_m$ and~$I^E_{s,t}$ is given by the
commutativity of the diagram
$$
\commdiag{
k& \isomorphism{I^E_{s,t}}& k\cr
\leftlabel{\gamma_m}\upar\rightlabel\sim& &
\leftlabel\sim\downar\rightlabel{\Norm_{E/k}v}\cr
k& \isomorphism{t_m}& \Norm_{E/k}\L(-D_m)\rlap.}
$$
This proves that the element of~$k^\times$ output of the last step is
indeed $I^E_{s,t}$.

It is straightforward to check that the running time of the algorithm,
measured in operations in~$k$, is polynomial in $\deg\L$
and~$m$.\endanalysis

\algorithm (Frey--R\"uck pairing). Let $X$ be a projective curve over
a finite field~$k$, let $n$ be an integer dividing $\#k^\times$, and
let $x$ and~$y$ be elements of~$J(k)$ with $nx=0$.  Given the
$k$-algebra~$S_X^{(7)}$ and subspaces $\gl(X,\L_X^{\otimes2}(-D))$
and~$\gl(X,\L_X^{\otimes2}(-E^-))$ of~$\gl(X,\L_X^{\otimes2})$
representing $x$ and~$y$, this algorithm outputs the element
$[x,y]_n\in\mu_n(k)$.

\label{frey-rueck-pairing}

\step Find an anti-addition chain $(a_0,a_1,\ldots,a_m)$ for~$n$.

\step Choose any non-zero global section $u$ of~$\L_X$, and let $D_0$
denote its divisor.  Compute the space
$$
\gl(X,\L_X^{\otimes2}(-D_0))=u\gl(X,\L_X).
$$
Write $D_1=D$.

\step Use Algorithm~\ref{addflip} to compute effective divisors $D_2$,
$D_3$, \dots, $D_m$ of degree~$\deg\L_X$, represented as the spaces
$\gl(X,\L_X^{\otimes2}(-D_l))$, and non-zero global sections $s_2$,
$s_3$, \dots, $s_m$ of~$\L_X^{\otimes3}$ such that the line bundle
$\L_X^{\otimes3}(-D_{i(l)}-D_{j(l)}-D_l)$ is trivial and
$$
\div(s_l)=D_{i(l)}+D_{j(l)}+D_l.
$$

\step Using Algorithm~\ref{zero-test}, verify that $\L_X(-D_m)$ is
trivial and find a non-zero global section $v$ of~$\L_X(-D_m)$.

\step Choose a non-zero global section~$w$ of~$\L_X$, let $E^+$ denote
its divisor, and compute
$$
\gl(X,\L_X^{\otimes2}(-E^+))=w\gl(X,\L_X).
$$

\step Compute $I^{E^+}_{s,t^+}$ and~$I^{E^-}_{s,t^-}$, viewed as
elements of~$k^\times$, using Algorithm~\ref{algorithm-Est}, where
$t^+$ and~$t^-$ are certain trivialisations chosen by that algorithm.

\step Output $(I^{E^+}_{s,t^+}/I^{E^-}_{s,t^-})^{\#k^\times/n}$.

\endalgorithm

\analysis The correctness of this algorithm follows from
Lemma~\ref{frey-rueck-via-norm}.  The running time is polynomial in
$\deg\L_X$, $\log\#k$ and~$\log n$.\endanalysis

\subsection{Finding relations between torsion points}

\label{finding-relations}

Let $X$ be a projective curve over a finite field~$k$, let $J$ be its
Jacobian, and let $l$ be a prime number different from the
characteristic of~$k$.  We will show how to find all the $\F_l$-linear
relations between given elements of~$J[l](k)$.  In particular, given a
basis $(b_1,\ldots,b_n)$ for a subspace~$V$ of~$J[l](k)$ and another
point~$x\in J[l](k)$, this allows us to check whether $x\in V$, and if
so, express $x$ as a linear combination of $(b_1,\ldots,b_n)$.

Let $k'$ be an extension of~$k$ containing a primitive $l$-th root of
unity.  It is well known that the problem just described can be
reduced, via the Frey--R\"uck pairing, to the discrete logarithm
problem in the group~$\mu_l(k')$.
Algorithm~\ref{algorithm-finding-relations} below makes this precise.
We begin with an estimate for the number of elements needed to
generate a finite-dimensional vector space over a finite field with
high probability.

\proclaim Lemma. Let $\F$ be a finite field, and let $V$ be an
$\F$-vector space of finite dimension $d$.  Let $\alpha$ be a real
number with $0<\alpha<1$, and write
$$
m=\cases{0& if\/ $d=0$;\cr
\displaystyle
d-1+\left\lceil{\log{1\over 1-\alpha^{1/d}}\over\log\#\F}\right\rceil&
if\/ $d>0$.}
$$
If $v_1$, \dots, $v_m$ are uniformly random elements of~$V$, the
probility that $V$ is generated by $v_1$, \dots, $v_m$ is at least
$\alpha$.

\label{vector-space-generation}

\proof Fix a basis of~$V$.  The matrix of the linear map
$$
\eqalign{
\F^m&\longrightarrow V\cr
(c_1,\ldots,c_m)&\mapsto\sum_{i=1}^m c_iv_i}
$$
is a uniformly random $d\times m$-matrix over~$\F$.  The probability
that it has rank~$d$ is the probability that its rows (which are
uniformly random elements of~$\F^m$) are linearly independent.  This
occurs with probability
$$
\eqalign{
p&={(\#\F^m-1)(\#\F^m-\#\F)\cdots(\#\F^m-\#\F^{d-1})\over\#\F^{dm}}\cr
&\ge{(\#\F^m-\#\F^{d-1})^d\over\#\F^{dm}}\cr
&=\bigl(1-(\#\F)^{-(m-d+1)}\bigr)^d}
$$
The choice of~$m$ implies that $p\ge\alpha$.\endproof

\remark The integer~$m$ defined in Lemma~\ref{vector-space-generation}
is approximately $d-1+{\log d\over\log\#F}$, in the sense that for any
fixed $\alpha$ the difference is bounded for $d\ge1$.

\algorithm (Relations between torsion points).  Let $X$ be a
projective curve over a finite field~$k$, let $J$ be its Jacobian, and
let $l$ be a prime number different from the characteristic of~$k$.
Let $x_1$, \dots, $x_n$ be elements of~$J[l](k)$.  Given the
$k$-algebra $S_X^{(h)}$ for some $h\ge7$ and subspaces
$\gl(X,\L_X^{\otimes2}(-D_i))$ of~$\gl(X,\L_X^{\otimes2})$
representing $x_i$ for $1\le i\le n$, this algorithm outputs an
$\F_l$-basis for the kernel of the natural map
$$
\eqalign{
\Sigma\colon\F_l^n&\longrightarrow J[l](k)\cr
(c_1,\ldots,c_n)&\longmapsto\sum_{i=1}^n c_ix_i.}
$$
The algorithm depends on a parameter $\alpha\in(0,1)$.

\step Generate a minimal extension $k'$ of~$k$ such that $k'$
contains a primitive $l$-th root of unity $\zeta$.  Let
$$
\lambda\colon \mu_l(k')\isom\F_l
$$
denote the corresponding discrete logarithm, i.e.\ the unique
isomorphism of one-dimensional $\F_l$-vector spaces sending $\zeta$
to~1.

\step Define an integer $m\ge0$ by
$$
m=\cases{
0& if $n=0$;\cr
\displaystyle
n-1+\left\lceil{\log{1\over 1-\alpha^{1/n}}\over\log l}\right\rceil&
if $n>0$.}
$$

\step Choose $m$ uniformly random elements $y_1$, \dots, $y_m$
in~$J(k')$ as described in~\S\ref{random-elements-picard-group}; their
images in~$J(k')/lJ(k')$ are again uniformly distributed.

\label{pick-y}

\step Compute the $m\times n$-matrix
$$
M=\left(\lambda([y_i,x_j]_l)\right)\quad (1\le i\le m,\ 1\le j\le n)
$$
with coefficients in~$\mu_l(k')$, where the
pairing~$[\blank\mathord,\blank]_l$ is evaluated using
Algorithm~\ref{frey-rueck-pairing} and the isomorphism~$\lambda$ is
evaluated using some algorithm for computing discrete logarithms
in~$\mu_l(k)$.

\step Compute an $\F_l$-basis $(b_1,\ldots,b_r)$ for the kernel
of~$M$.

\step If $\Sigma(b_1)=\ldots=\Sigma(b_r)=0$, output
$(b_1,\ldots,b_r)$ and stop.

\step Go to step \ref{pick-y}.

\endalgorithm

\label{algorithm-finding-relations}

\analysis We write $V$ for the image of~$\Sigma$ and $V'$ for the
quotient of~$J(k')/lJ(k')$ by the annihilator of~$V$ under the
pairing~$[\blank\mathord,\blank]_l$.  Then we have an induced
isomorphism
$$
V\isom\Hom_{\F_l}(V',\mu_l(k')).
$$
Consider the map
$$
\eqalign{
\Sigma'\colon\F_l^m&\longrightarrow V'\cr
(c_1,\ldots,c_m)&\longmapsto\sum_{i=1}^m c_iy_i.}
$$
Now we have a commutative diagram
$$
\commdiag{
\F_l^n& \longrightarrow& \Hom_{\F_l}(\F_l^m,\mu_l(k'))\cr
\leftlabel{\Sigma}\downar& & \upar\rightlabel{f\mapsto f\circ\Sigma'}\cr
V& \isom& \Hom_{\F_l}(V',\mu_l(k'))}
$$
We identify $\mu_l(k')$ with~$\F_l$ using the isomorphism~$\lambda$
and equip $\Hom_{\F_l}(\F_l^m,\mu_l(k'))$ with the dual basis of the
standard basis of~$\F_l^m$.  Then the top arrow in the diagram is
given by the matrix~$M$ defined in step~4.  This means that we have an
inclusion
$$
\ker\Sigma\subseteq\ker M.
$$
In step~6 we check whether this inclusion is an equality.  The
surjectivity of~$\Sigma$ implies that this is the case and only if the
rightmost map in the diagram is injective, i.e.\ if and only if
$\Sigma'$ is surjective.  Since $\dim_{\F_l} V\le n$, this happens
with probability at least $\alpha$ by
Lemma~\ref{vector-space-generation}.  Therefore steps 3--7 are
executed at most $1/\alpha$ times on average.  This implies that (for
fixed $\alpha$) the algorithm runs in time polynomial in $\genus_X$,
$\log\#k$, $l$ and~$n$.\endanalysis

\remarks (1)\enspace If we know an upper bound for the dimension of
the $\F_l$-vector space generated by the $x_i$, then we can use this
upper bound instead of~$n$ in the expression for~$m$ in step~2.
\smallskip\noindent
(2)\enspace It does not matter much what algorithm we use for
computing the discrete logarithm in~$\mu_l(k')$, since the running
time of Algorithm~\ref{algorithm-finding-relations} is already
polynomial in~$l$.  For example, we can simply tabulate the
function~$\lambda$.

\subsection{The Kummer map on a divisible group}

\label{kummer-map}

Let $k$ be a finite field of cardinality~$q$, and let $l$ be a prime
number.  Let $\G$ be an \'etale $l$-divisible group over~$k$.
(The \'etaleness is automatic if $l$ is different from the
characteristic of~$k$.)  We denote by $\Frob_q\colon\G\to\G$ the
($q$-power) Frobenius endomorphism of~$\G$; this is an automorphism
because of the assumption that $\G$ is \'etale.

For any non-negative integer~$n$ such that all the points of~$\G[l^n]$
are $k$-rational, the {\it Kummer map\/} of order~$l^n$ on~$\G$
over~$k$ is the isomorphism
$$
\eqalign{
K_{l^n}^{\G/k}\colon\G(k)/l^n\G(k)&\isom\G[l^n](k)\cr
x&\longmapsto\Frob_q(y)-y,}
$$
where $y$ is any point of~$\G$ over an algebraic closure of~$k$ such
that $l^ny$ is a lift of~$x$ to $\G(k)$.

Let $\chi\in\Z_l[t]$ be the characteristic polynomial of the Frobenius
automorphism of~$\G$ on (the Tate module of)~$\G$.  Then the element
$t\bmod\chi$ of~$\Z_l[t]/(\chi)$ is invertible.  Let $n$ be any
non-negative integer, and let $a$ be a positive integer such that
$$
t^a=1\quad\hbox{in }(\Z_l[t]/(l^n,\chi))^\times.
$$
Then $t^a-1$ is divisible by~$l^n$ in~$\Z_l[t]/(\chi)$, and we let
$h_a$ be the unique element of~$\Z_l[t]/(\chi)$ such that
$$
t^a-1=l^nh_a\in\Z_l[t]/(\chi).
$$
By the Cayley--Hamilton theorem, $\Z_l[t]/(\chi)$ acts on~$\G$ with
$t$ acting as~$\Frob_q$.  The above identity therefore implies that
$$
\Frob_q^a-1=l^nh_a(\Frob_q)\quad\hbox{on }\G.
$$
Let $k_a$ be an extension of~$k$ with
$$
[k_a:k]=a.
$$
Then $\G[l^n]$ is defined over~$k_a$, and we can express the Kummer map
over~$k_a$ in terms of the Frobenius endomorphism over~$k$ as
$$
\eqalign{
K_{l^n}^{\G/k_a}\colon\G(k_a)/l^n\G(k_a)&\isom\G[l^n](k_a)\cr
x&\longmapsto h_a(\Frob_q)(x).}
$$
In~\S\ref{computing-with-torsion-points} we are going to apply this
to a certain $l$-divisible subgroup of the $l$-power torsion of the
Jacobian of a projective curve over~$k$.

\subsection{Computing the $l$-torsion in the Picard group}

\label{computing-with-torsion-points}

Let $X$ be a projective curve over~$k$, and let $J$ be its Jacobian.
Let $\Frob_q$ denote the Frobenius endomorphism of~$J$ over~$k$, and
let $\chi\in\Z[t]$ be the characteristic polynomial of~$\Frob_q$.

Let $l$ be a prime number different from the characteristic of~$k$.
We are going to apply the results of~\S\ref{kummer-map} to a certain
$l$-divisible subgroup~$\G$ of the group~$J[l^\infty]$ of $l$-power
torsion points of~$J$.  This $\G$ is defined as follows.  Let $\bar
f=(t-1)^b$ be the largest power of~$t-1$ dividing $\chi\bmod l$, so
that $\chi\bmod l$ has the factorisation
$$
(\chi\bmod l)=\bar f\cdot\bar f^\perp
$$
in coprime monic polynomials in~$\F_l[t]$.  Hensel's lemma implies
that this factorisation can be lifted uniquely to a factorisation
$$
\chi=f\cdot f^\perp,
$$
where $f$ and~$f^\perp$ are coprime monic polynomials in~$\Z_l[t]$.
The Chinese remainder theorem gives a decomposition
$$
\Z_l[t]/(\chi)\isom\Z_l[t]/(f)\times\Z_l[t]/(f^\perp),
\eqnumber{chinese-decomposition}
$$
which in turn induces a decomposition
$$
J[l^\infty]\cong\G\times\G^\perp
$$
of $l$-divisible groups.  We note that $\G$ is of rank~$b$ and that
$f$ is the characteristic polynomial of~$\Frob_q$ on~$\G$.  Let $a$ be
a positive integer such that
$$
t^a=1\quad\hbox{in }(\F_l[t]/\bar f)^\times,
\eqnumber{a-multiple-order}
$$
let $h_a$ be the unique element of~$\Z_l[t]/(f)$ such that
$$
t^a-1=lh_a\in\Z_l[t]/(f),
\eqnumber{definition-h_a}
$$
and let $k_a$ be an extension of degree~$a$ of~$k$.  All the points
of~$\G[l]$ are $k_a$-rational, and the $b$-dimensional $\F_l$-vector
space $\G[l](k_a)$ is the generalised eigenspace corresponding to the
eigenvalue~1 of~$\Frob_q$ inside the $\F_l$-vector space of points
of~$J[l]$ over an algebraic closure of~$k_a$.  In particular, we have
the identity
$$
J[l](k)=\{x\in\G[l](k_a)\mid\Frob_q(x)=x\}.
$$
As explained in \S\ref{kummer-map}, the map
$$
\eqalign{
\G(k_a)/l\G(k_a)&\isom\G[l](k_a)\cr
x&\longmapsto h_a(\Frob_q)(x)}
$$
is well-defined and equal to the Kummer isomorphism
$$
K_l^{\G/k_a}\colon\G(k_a)/l\G(k_a)\isom\G[l](k_a)
$$
of order~$l$.

The above results give us a way of generating uniformly random
elements of the $\F_l$-vector space~$\G[l](k_a)$.  We factor
$\#J(k_a)$ as
$$
\#J(k_a)=l^{c_a}m_a
$$
with $c_a\ge0$, $m_a\ge 1$ and $l\nmid m_a$.  Let $e$ be the
idempotent in~$\Z_l[t]/(\chi)$ corresponding to the element $(1,0)$ on
the right-hand side of~\eqref{chinese-decomposition}.  Composing the
maps
$$
J(k_a)\morphism{m_a}J[l^\infty](k_a)
\morphism{e(\Frob_q)}\G(k_a)
\longrightarrow\G(k_a)/l\G(k_a)
\morphism{h_a(\Frob_q)}\G[l](k_a)
\eqnumber{random-G[l]}
$$
we get a surjective group homomorphism from~$J(k_a)$ to~$\G[l](k_a)$.
We can use this map to convert uniformly random elements of~$J(k_a)$
into uniformly random elements of~$\G[l](k_a)$, provided we know $e$
and~$h_a$ to sufficient $l$-adic precision.  It is clear that to
compute the Kummer map we only need to know the image of~$h_a$
in~$\Z_l[t]/(f,l)=\F_l[t]/((t-1)^b)$.  Since $\G(k_a)$ can be
identified with a subgroup of~$\#J(k_a)$, it is annihilated
by~$l^{c_a}$, and we have
$$
J[l^\infty](k_a)=J[l^{c_a}](k_a)\quad\hbox{and}\quad
\G(k_a)=\G[l^{c_a}](k_a).
$$
This implies that it suffices to know $e$ to precision~$O(l^{c_a})$.

Let us check that there is a reasonably small $a$ for
which \eqref{a-multiple-order} holds.  For any non-negative
integer~$\gamma$ the identity
$$
t^{l^\gamma}-1=(t-1)^{l^\gamma}
$$
holds in~$\F_l[t]$, and the right-hand side maps to zero in
$\F_l[t]/(t-1)^b$ if and only if $l^\gamma\ge b$.  Since $l$ is a
prime number, we conclude that the order of~$t$ in~$\F_l[t]/((t-1)^b)$
equals $l^\gamma$, where $\gamma$ is the least non-negative integer
such that $l^\gamma\ge b$.

\algorithm (Computing the $l$-torsion of the Picard group). Let $X$ be
a projective curve over a finite field~$k$ with $q$ elements, let $J$
be its Jacobian, and let $l$ be a prime number different from the
characteristic of~$k$.  Given the $k$-algebra $S_X^{(7)}$ and the
characteristic polynomial $\chi$ of the Frobenius endomorphism of~$J$
over~$k$, this algorithm outputs an $\F_l$-basis for~$J[l](k)=(\Pic
X)[l]$.  The algorithm depends on a parameter~$\alpha\in(0,1)$.

\step Factor $\chi\bmod l$ in~$\F_l[t]$ as
$$
(\chi\bmod l)=\bar f\cdot\bar f^\perp,
$$
where $\bar f$ is the greatest power of $t-1$ dividing $\chi\bmod l$,
say $\bar f=(t-1)^b$, and lift this to a factorisation
$$
\chi=f\cdot f^\perp
$$
in coprime monic polynomials in~$\Z_l[t]$.

\step Compute the non-negative integer $r$ defined by
$$
r=\cases{0& if $b=0$;\cr
\displaystyle
b-1+\left\lceil{\log{1\over1-\alpha^{1/b}}\over\log l}\right\rceil&
if $b\ge 1$.}
$$

\step Define $a=l^\gamma$, where $\gamma$ is the least non-negative
integer such that $l^\gamma\ge b$.  Generate a finite extension $k_a$
of degree~$a$ of~$k$.  Factor $\#J(k_a)$ as
$$
\#J(k_a)=l^{c_a}m_a\quad\hbox{with }l\nmid m_a.
$$
Compute the image of the idempotent~$e$ in~$(\Z/l^{c_a}\Z)[t]/(\chi)$
using the extended Euclidean algorithm, and compute the image of~$h_a$
in~$\F_l[t]/((t-1)^b)$ using the definition~\eqref{definition-h_a}
of~$h_a$.

\step Generate $r$ uniformly random elements of~$J(k_a)$ as explained
in~\S\ref{random-elements-picard-group}, and map them to elements
$x_1,\ldots, x_r\in\G[l](k_a)$ using the
homomorphism~\eqref{random-G[l]}.

\label{pick-x_i}

\step Using Algorithm~\ref{finding-relations}, compute a basis for the
kernel of the $\F_l$-linear map
$$
\eqalign{
\Sigma\colon\F_l^r&\longrightarrow \G[l](k_a)\cr
(c_1,\ldots,c_r)&\longmapsto\sum_{i=1}^r c_ix_i.}
$$
If the dimension of this kernel is greater than $r-b$, go to
step~\ref{pick-x_i}.

\step Use the $\F_l$-linear relations between $x_1$, \dots, $x_r$
computed in the previous step to find a subsequence $(y_1,\ldots,y_b)$
of~$(x_1,\ldots,x_r)$ that is an $\F_l$-basis of~$\G[l](k_a)$.

\step Let $M$ be the matrix with respect to the basis
$(y_1,\ldots,y_b)$ of the $\F_l$-linear automorphism of~$\G[l](k_a)$
induced by the Frobenius endomorphism~$\Frob_q$ of~$J$ over~$k$.
Compute $M$ by computing $\Frob_q(y_i)$ for $i=1$, \dots, $b$ using
Algorithm~\ref{compute-frobenius-map} and then applying
Algorithm~\ref{finding-relations} to express the $\Frob_q(y_i)$ as
linear combinations of the $y_i$.

\step Compute a basis for the kernel of $M-I$, where $I$ is the
$b\times b$ identity matrix.  Map the basis elements to elements
$z_1$, \dots, $z_t$ of~$\G[l](k_a)$ using the injective homomorphism
$$
\eqalign{
\F_l^b&\longrightarrow\G[l](k_a)\cr
(a_1,\ldots,a_b)&\longmapsto\sum_{i=1}^b a_iy_i.}
$$
Output $(z_1,\ldots,z_t)$.

\endalgorithm

\label{compute-l-torsion}

\analysis The definition of~$a$ implies that $a$ equals the order
of~$t$ in~$(\F_l[t]/(t-1)^b)^\times$, and $J[l](k)$ equals the kernel
of $\Frob_q-\id$ on~$\G[l](k_a)$, as remarked before.  The elements
$x_1$, \dots, $x_r$ of~$\G[l](k_a)$ are uniformly random by the fact
that \eqref{random-G[l]} is a homomorphism.  By
Lemma~\ref{vector-space-generation}, they generate the $b$-dimensional
$\F_l$-vector space $\G[l](k_a)$ with probability at least $\alpha$.
The definition of~$a$ also implies that
$$
a\le\max\{1,2\genus_Xl-1\},
$$
while the ``class number formula''~\eqref{class-number-formula} gives
the upper bound
$$
\eqalign{
c_a&\le{\log\#J(k_a)\over\log l}\cr
&\le{2\genus_X\log\bigl(1+q^{a/2}\bigr)\over\log l}.}
$$
This shows that $c_a$ is bounded by a polynomial in $\genus_X$, $\log
q$ and~$l$.  For fixed $\alpha$ we therefore reach step~6 in expected
polynomial time in $\deg\L_X$, $\log q$ and~$l$.  In steps 6--8 we
compute a basis for the kernel of~$\Frob_q-\id$, which is $J[l](k)$.
We conclude that the algorithm is correct and runs in probabilistic
polynomial time in $\deg\L_X$, $\log q$ and~$l$.\endanalysis

\remark The elements $z_j$ output by the preceding algorithm are
defined over~$k$.  In general, it seems unclear how to generate
$k$-vector spaces (instead of $k_a$-vector spaces) representing them.
However, if we know a $k$-rational point on~$X$, then we can use
Algorithm~\ref{descent} to accomplish this.

\unnumberedsection{References}

\parindent=\normalparindent
\advance\parskip by1ex

\reference{Adleman-Lenstra} L. M. {\sc Adleman} and H. W. {\sc
Lenstra}, Jr., Finding irreducible polynomials over finite fields.
In: {\sl Proceedings of the Eighteenth Annual ACM Symposium on Theory
of Computing (Berkeley, CA, 1986)\/}, 350--355.  Association for
Computing Machinery, New York, 1986.

\reference{Bosman} J. G. {\sc Bosman}, {\sl Explicit computations with
modular Galois representations\/}.  \PhD thesis, Universiteit Leiden,
2008.

\reference{thesis} P. J. {\sc Bruin}, {\sl An algorithm for computing
modular Galois representations\/}.  \PhD thesis, Universiteit Leiden,
2010, in preparation.

\reference{Couveignes: Linearizing torsion classes} J.-M. {\sc
Couveignes}, Linearizing torsion classes in the Picard group of
algebraic curves over finite fields.  {\it Journal of Algebra\/} {\bf
321} (2009), 2085--2118.

\reference{Diem} C. {\sc Diem}, {\sl On arithmetic and the discrete
logarithm problem in class groups of curves\/}.
Habi\-li\-ta\-tions\-schrift, Universit\"at Leipzig, 2008.

\reference{Diem: article} C. {\sc Diem}, On the discrete logarithm
problem in class groups of curves, to appear.

\reference{Eberly-Giesbrecht} W. {\sc Eberly} and M. {\sc Giesbrecht},
Efficient decomposition of associative algebras over finite fields.
{\it Journal of Symbolic Computation\/} {\bf 29} (2000), 441--458.

\reference{book} J.-M. {\sc Couveignes} and S. J. {\sc Edixhoven}
(editors), {\sl Computational aspects of modular forms and Galois
representations\/}.  Princeton University Press, to appear.

\reference{compcoefs} S. J. {\sc Edixhoven} (with J.-M. {\sc
Couveignes}, R. S. {\sc de Jong}, F. {\sc Merkl} and J. G. {\sc
Bosman}), On the computation of coefficients of a modular form.
Preprint, 2006/2009.\hfill\break
Available online: {\tt http://arxiv.org/abs/math.NT/0605244}.

\reference{Frey-Rueck} G. {\sc Frey} and H.-G. {\sc R\"uck}, A remark
concerning $m$-divisibility and the discrete logarithm in the divisor
class group of curves.  {\it Mathematics of Computation\/} {\bf 62}
(1994), 865--874.

\reference{Hartshorne} R. {\sc Hartshorne}, {\sl Algebraic
Geometry\/}.  Springer-Verlag, New York, 1977.

\reference{Khuri-Makdisi: Linear algebra algorithms} K. {\sc
Khuri-Makdisi}, Linear algebra algorithms for divisors on an algebraic
curve.  {\it Mathematics of Computation\/} {\bf 73} (2004), no.~245,
333--357.\hfill\break
Available online: {\tt http://arxiv.org/abs/math.NT/0105182}.

\reference{Khuri-Makdisi} K. {\sc Khuri-Makdisi}, Asymptotically fast
group operations on Jacobians of general curves.  {\it Mathematics of
  Computation\/} {\bf 76} (2007), no.~260, 2213--2239.\hfill\break
Available online: {\tt http://arxiv.org/abs/math.NT/0409209}.

\reference{Lazarsfeld} R. {\sc Lazarsfeld}, A sampling of vector
bundle techniques in the study of linear series.  In: M. {\sc
Cornalba}, X. {\sc Gomez-Mont} and A. {\sc Verjovsky} (editors), {\sl
Lectures on Riemann Surfaces (Trieste, 1987)\/}, 500--559.  World
Scientific Publishing, Teaneck, NJ, 1989.

\reference{Rabin} M. O. {\sc Rabin}, Probabilistic algorithms in
finite fields.  {\it SIAM Journal on Computing\/} {\bf 9} (1980),
no.~2, 273--280.

\reference{Schaefer} E. F. {\sc Schaefer}, A new proof for the
non-degeneracy of the Frey--R\"uck pairing and a connection to
isogenies over the base field.  In: T. {\sc Shaska} (editor), {\sl
Computational Aspects of Algebraic Curves\/} (Conference held at the
University of Idaho, 2005), 1--12.  Lecture Notes Series in Computing
{\bf 13}.  World Scientific Publishing, Hackensack, NJ, 2005.

\reference{SGA 4 III} {\sl Th\'eorie des topos et cohomologie \'etale
  des sch\'emas} (SGA 4).  Tome 3 (expos\'es IX \`a XIX)\null.
S\'eminaire de G\'eom\'etrie Alg\'ebrique du Bois-Marie 1963--1964,
dirig\'e par M. {\sc Artin}, A. {\sc Gro\-then\-dieck} et J.-L. {\sc
  Verdier}, avec la collaboration de P. {\sc Deligne} et B. {\sc
  Saint-Donat}.  Lecture Notes in Mathematics {\bf 305}.
Springer-Verlag, Berlin/Heidelberg/New York, 1973.

\reference{Stein} W. A. {\sc Stein}, {\sl Modular Forms, a
Computational Approach\/}.  With an appendix by P. E. {\sc Gunnells}.
American Mathematical Society, Providence, RI, 2007.

\vskip1cm
\vbox{
\hbox{Peter Bruin}
\hbox{Universiteit Leiden}
\hbox{Mathematisch Instituut}
\hbox{Postbus 9512}
\hbox{2300 RA \ Leiden}
\hbox{Netherlands}
\hbox{\tt pbruin@math.leidenuniv.nl}
}

\bye